\documentclass[11pt]{article}
\usepackage{amsfonts}
\usepackage{bbm}
\usepackage[top=1in, bottom=1in, left=1in, right=1in]{geometry}
\usepackage{amsmath,booktabs,ctable,threeparttable}
\usepackage{amssymb,amsfonts,boxedminipage}
\usepackage{amsthm}
\usepackage{algorithm}
\usepackage{algorithmic}
\usepackage{epsfig,graphicx,picins,picinpar,subfigure,epstopdf,float}
\usepackage{multirow}
\usepackage{verbatim}
\usepackage{bm}
\usepackage[colorlinks,
            linkcolor=red,
            anchorcolor=blue,
            citecolor=blue
            ]{hyperref}

\numberwithin{equation}{section}
\newtheorem{thm}{Theorem}[section]

\theoremstyle{definition}

\begin{document}

\title{On Regularizing Effects of MINRES and MR-II
for Large Scale Symmetric Discrete Ill-posed Problems\thanks
{This work was supported in part by the National Science Foundation of China
(No. 11371219).}}

\author{\small Yi Huang$^a$, Zhongxiao Jia$^b$\\
        \small$^a$ Department of Mathematical Sciences, Tsinghua
        University, 100084 Beijing, China\\
       \small $^b$ Department of Mathematical Sciences, Tsinghua
       University, 100084 Beijing, China,\\
        \small  E-mail:\ \ huangyi10@mails.tsinghua.edu.cn,\ \
        jiazx@tsinghua.edu.cn\\}
\date{}
\maketitle

\begin{abstract}
For large scale symmetric discrete ill-posed problems, MINRES and MR-II are
often used iterative regularization solvers. We call a regularized solution
best possible if it is at least as accurate as the best regularized solution
obtained by the truncated singular value decomposition (TSVD) method.
In this paper, we analyze
their regularizing effects and establish the following results:
(i) the filtered SVD expression are derived for the regularized solutions by MINRES;
(ii) a hybrid MINRES that uses explicit regularization
within projected problems is needed to compute a best possible
regularized solution to a given ill-posed problem;
(iii) the $k$th iterate by MINRES is more accurate
than the $(k-1)$th iterate by MR-II until the semi-convergence
of MINRES, but MR-II has globally better regularizing effects than MINRES;
(iv) bounds are obtained for the
2-norm distance between an underlying $k$-dimensional Krylov subspace and
the $k$-dimensional dominant eigenspace. They show that MR-II has better
regularizing effects for severely and moderately ill-posed
problems than for mildly ill-posed problems, and
a hybrid MR-II is needed to get a best possible regularized solution
for mildly ill-posed problems; (v) bounds are derived for the entries generated
by the symmetric Lanczos process that MR-II is based on,
showing how fast they decay. Numerical experiments confirm our assertions.
Stronger than our theory, the regularizing effects of MR-II are experimentally
shown to be good enough to obtain best possible regularized
solutions for severely and moderately ill-posed problems.

\smallskip

\end{abstract}

\section{Introduction}
\label{intro}
Consider the large scale discrete linear ill-posed problem
\begin{equation}\label{eq1}
    Ax=b,\ \ \ A\in \mathbb{R}^{n\times n},
  \ b\in \mathbb{R}^{n},
\end{equation}
where $A$ is symmetric and extremely ill conditioned with
its singular values decaying gradually to zero without a noticeable gap.
This kind of problem arises in many science and engineering areas \cite{hansen10}.
The right-hand side $b$ is noisy and typically affected
by a white noise, caused by measurement, modeling or discretization
errors, i.e.,
\begin{equation*}\label{}
  b=\hat{b}+e,
\end{equation*}
where $e\in \mathbb{R}^{n}$ represents a  white noise vector and
$\hat{b}\in \mathbb{R}^{n}$ denotes the noise-free right-hand side, and it
is supposed that $\|e\|<\|\hat{b}\|$.
Because of the presence of noise $e$ in $b$ and the high ill-conditioning of
$A$, the naive solution $x_{naive}=A^{-1}b$ of \eqref{eq1} is
far from the true solution $x_{true}=A^{-1}\hat{b}$ and meaningless.
Therefore, one needs to use regularization to determine a
regularized solution so that it is close to $x_{true}=A^{-1}\hat{b}$
as much as possible \cite{hansen98,hansen10}.

For $A$ symmetric, its SVD is closely related to its spectral decomposition
as follows:
\begin{equation}\label{decom}
  A=V\Lambda V^T=V\Omega\Sigma V^T=U\Sigma V^T,
\end{equation}
where $U=(u_1,u_2,\ldots,u_n)=V\Omega$ and $V=(v_1, v_2, \ldots, v_n)$ are
orthogonal, whose columns are the left and right singular vectors of $A$,
respectively, the diagonal matrix $\Sigma =
\mathrm{diag}(\sigma_1,\sigma_2,\ldots,\sigma_n)$
with the singular values labeled as $\sigma_1> \sigma_2>\cdots > \sigma_n>0$,
$\Omega=\mathrm{diag}(\pm1)$ is a signature matrix such that
$\sigma_i=|\lambda_i|$ with the $\lambda_i$ the eigenvalues of $A$, and
$\Lambda=\mathrm{diag}(\lambda_1,\lambda_2,\ldots,\lambda_n)$. Obviously,
$u_i=\pm v_i$ with the $\pm$ sign depending on $\Omega$.
With \eqref{decom}, we can express the naive solution of \eqref{eq1}
as
\begin{equation}\label{eq4}
  x_{naive}=\sum\limits_{i=1}^{n}\frac{v_i^{T}b}{\lambda_i}v_i =
  \sum\limits_{i=1}^{n}\frac{v_i^{T}\hat{b}}{\lambda_i}v_i +
  \sum\limits_{i=1}^{n}\frac{v_i^{T}e}{\lambda_i}v_i
  =x_{true}+\sum\limits_{i=1}^{n}\frac{v_i^{T}e}{\lambda_i}v_i.
\end{equation}

Throughout the paper, we assume that $\hat{b}$
satisfies the discrete Picard condition \cite{hansen98,hansen10}: On average, the
coefficients $|u_i^{T}\hat{b}|=|v_i^{T}\hat{b}|$ decay faster than the
singular values $\sigma_i=|\lambda_i|$.
This is a necessary hypothesis
that controls the size of $\|x_{true}\|$ and makes regularization
possible to find meaningful approximations to $x_{true}$
\cite{hansen98,hansen10}.
For the sake of simplicity, we assume that they satisfy a widely used model
\cite{hansen98,hansen10}:
\begin{equation}\label{picard}
  \mid v_i^T \hat b\mid=|\lambda_i|^{1+\beta},\ \ \beta>0,\ i=1,2,\ldots,n.
\end{equation}

Similar to the truncated
SVD (TSVD) method \cite{hansen98,hansen10}, for $A$ symmetric,
a truncated spectral decomposition method obtains the TSVD
regularized solutions
\begin{equation}\label{tsvdsolution}
x_{k}=\sum\limits_{i=1}^{k}\frac{v_i^{T}b}{\lambda_i}v_i=
\sum\limits_{i=1}^{k}\frac{u_i^{T}b}{\sigma_i}v_i=A_{k}^{\dagger}b,\
k=1,2,\ldots,n,
\end{equation}
where $A_{k}=U_{k}\Sigma_{k} V_{k}^T$
with $U_{k}$ and $V_{k}$ the first $k$ columns of $U$ and $V$, respectively,
$\Sigma_{k}={\rm diag}(\sigma_1,\ldots,
\sigma_{k})$, and $\dagger$ the Moore-Penrose generalized inverse of a matrix.
Obviously, $x_{k}$ is the minimum 2-norm least squares solution of
the perturbed problem that replaces $A$ in \eqref{eq1} by its best rank $k$
approximation $A_{k}$.

Let $k_0$ denote the transition point such that $\mid v_{k_0}^T \hat b\mid>
\mid v_{k_0+1}^T e\mid$ and $\mid v_{k_0+1}^T \hat b\mid\approx
\mid v_{k_0+1}^T e\mid$, which divides the eigenvalues $\lambda_i$ or
equivalently the singular values
$\sigma_i=|\lambda_i|$ into the dominant or large ones for $i\leq k_0$
and the small ones for $i>k_0$. It is known from \cite[p. 176]{hansen98} and
\cite[p. 86-88]{hansen10} that the best TSVD regularized
solution is $x_{k_0}$, which consists of the $k_0$ dominant SVD components of
$A$, i.e., the dominant spectral components corresponding to the first $k_0$
large eigenvalues in magnitude when $A$ is symmetric. A number of approaches
have been proposed for determining $k_0$, such as discrepancy principle,
the discrete L-curve criterion and the generalized cross validation (GCV);
see, e.g., \cite{bauer11,gazzola14,hansen98,kilmer03,reichel13}
for comparisons of the classical and new ones. In our numerical experiments, we
do this using the L-curve criterion in the TSVD method and Krylov iterative methods.

The TSVD method is important in its own right and plays a central role in
analyzing the standard-form Tikhonov regularization \cite{hansen98,hansen10}.
It and the standard-form Tikhonov regularization
expand their solutions in the basis vectors $v_i$
and produce very similar solutions with essentially the minimum
2-norm error; see \cite[p. 109-111]{hansen98} and
\cite[Sections 4.2 and 4.4]{hansen10}. Therefore,
the TSVD method can get a best regularized solution to \eqref{eq1}, and it
has long been used as a general-purpose reliable and efficient numerical
method for solving a small and/or moderate sized \eqref{eq1} \cite{hansen98,hansen10}.
As a result, we will take the TSVD
solution $x_{k_0}$ as a standard reference when assessing the regularizing
effects of iterative solvers and accuracy of iterates
under consideration in this paper.

For \eqref{eq1} is large, it is generally impractical to compute the spectral
decomposition of $A$. In this case, one typically solves it iteratively
via some Krylov subspace methods \cite{hansen98,hansen10}.
For \eqref{eq1} with a general matrix $A$,
the mathematically equivalent LSQR \cite{paige82} and CGLS \cite{bjorck96}
have been most commonly used for years and have been shown to have intrinsic
regularizing effects \cite{hanke01,hansen98,hansen10}.
They exhibit the semi-convergence; see \cite{hanke01}, \cite[p. 135]{hansen98},
\cite[p. 110]{hansen10}: The iterates
tend to be better and better approximations to the exact solution $x_{true}$
and their norms increase slowly and the residual norms decrease.
In later stages, however, the noise $e$ starts to deteriorate the iterates,
so that they will start to diverge from
$x_{true}$ and instead converge to $x_{naive}$,
while their norms increase considerably and the residual norms stabilize.
For LSQR, the semi-convergence is due to the fact that the projected
problem at some iteration starts to inherit the ill-conditioning
of \eqref{eq1}, that is, the noise progressively enters
the solution subspace, so that a small singular value of the projected
problem appears and the regularized solution is
deteriorated \cite{hansen98,hansen10}.

As far as an iterative solver for solving \eqref{eq1} is concerned,
a central problem is whether or not it has already obtained
a best possible regularized solution at semi-convergence.
Here, as defined in the abstract,
a best possible regularized solution means that it is at least as
accurate as the best TSVD solultion $x_{k_0}$. This problem has been
intensively studied but has had no definitive solutions. For Krylov iterative
solvers, their regularizing effects critically rely on how well the
underlying $k$-dimensional Krylov subspace captures the $k$ dominant
right singular vectors of $A$ \cite{hansen98,hansen10}.
The richer information the Krylov subspace contains on the $k$ dominant
right singular vectors, the less possible it is that
the resulting projected problem has a small singular value.
That is, the solvers capture the large SVD components of $A$
more effectively, and thus have better regularizing effects.

To precisely measure the regularizing effects, we introduce
the term of {\em full} or {\em partial} regularization.
If a pure iterative solver itself computes a best possible regularized
solution at semi-convergence, it is said
to have the full regularization; in this case, no additional regularization
is necessary. Otherwise, it is said to have the partial regularization; in
this case, a sophisticated hybrid variant is needed that combines the solver
with some additional regularization in order to improve the accuracy of
the regularized solution by the iterative solver at
semi-convergence \cite{hansen98,hansen10}. It appears that the regularizing
effects are closely related to the degree of ill-posedness of the problem.
To this end, we introduce the following definition of the degree of
ill-posedness, which follows Hofmann's book \cite{hofmann86} and has been
commonly used in the literature, e.g.,
\cite{hansen98,hansen10}: If there exists a positive real number
$\alpha$ such that the singular values satisfy $\sigma_j=
\mathcal{O}(j^{-\alpha})$, the problem is
termed as mildly or moderately ill-posed if $\alpha\le1$ or $\alpha>1$;
if $\sigma_j=\mathcal{O}(e^{-\alpha j})$ with $\alpha>0$ considerably,
$j=1,2,\ldots,n$, the problem is termed severely ill-posed.
Clearly, the singular values $\sigma_j$ of a
severely ill-posed problem decay exponentially at the same rate $e^{-\alpha}$,
while those of a moderately or mildly ill-posed problem decay more and more
slowly at the decreasing rate $\left(\frac{j}{j+1}\right)^{\alpha}$
approaching one with increasing $j$, which, for the same $j$, is
smaller for the moderately ill-posed problem than it
for the mildly ill-posed problem.

For $A$ symmetric, its $k$-dimensional dominant eigenspace is identical to
the $k$-dimensional dominant left and right singular subspaces. In this case,
MINRES and its variant MR-II are natural alternatives to
LSQR and CGLS \cite{hanke95,hanke96,hansen10}.
MR-II was originally designed for solving singular and inconsistent
linear systems, and it uses the starting vector
$Ab$ and restricts the resulting Krylov subspace
to the range of $A$. Thus, the iterates are orthogonal to the null
space of $A$, and MR-II computes
the minimum 2-norm least squares solution \cite{fischer96,hanke95}.
For \eqref{eq1}, we are not interested in such solution but a regularized
solution that is close to $x_{true}$ as much as possible. MINRES and MR-II
have been shown to have regularizing effects and exhibit
semi-convergence \cite{hansen08,jensen07,kilmer99}, and
MR-II usually provides better regularized solutions than MINRES. Intuitively,
this is because the noise $e$
in the initial Krylov vector $Ab$ is filtered by multiplication
with $A$ \cite{hanke01,jensen07}. Different implementations
associated with MR-II have been studied \cite{hanke95,neuman12}.

In this paper, we first prove that the MINRES iterates are filtered SVD
solutions, similar to the form of the LSQR iterates. Based on this result,
we show why MINRES, in general, has only the
partial regularization, independent of the degree of ill-posedness of \eqref{eq1}.
As a result, a hybrid MINRES that combines MINRES with a regularization method
applied to the lower dimensional projected problems should be used to
compute a best possible regularized solution; see \cite[Section 6.4]{hansen10}
for details.  Afterwards, we
take a closer look at the regularization of MINRES and MR-II in more
detail, which, from a new perspective, shows that MINRES
has only the partial regularization. We prove that, though MR-II has globally
better regularizing effects than MINRES, the $k$th MINRES
iterate is always more accurate than the $(k-1)$th MR-II iterate
until the semi-convergence of MINRES.
In a manner different from those used in \cite{jensen07,kilmer99},
we then analyze the regularizing effects of MR-II
and draw some definitive conclusions. We establish bounds
for the 2-norm distance between the underlying $k$-dimensional
Krylov subspace and the $k$-dimensional dominant eigenspace.
The bounds indicate that
the Krylov subspace better captures the $k$-dimensional
dominant eigenspace for severely and moderately ill-posed problems
than for mildly ill-posed problems. As a consequence, MR-II
has better regularizing effects for the first two kinds of problems
than for the third kind, for which MR-II has only the partial
regularization.
We then use the results to derive an estimate for the accuracy of
the rank $k$ approximation generated by MR-II to $A$.
Finally, we derive estimates for the entries generated
by the symmetric Lanczos process that
MR-II is based on, and show how fast they decay.

The paper is organized as follows. In Section \ref{SectionRev}, we describe
MINRES and MR-II. In Section \ref{SectionMin}, we prove that the MINRES
iterates are filtered SVD solutions, followed by an analysis on the regularizing
effects of MINRES. In Section \ref{SectionCom}, we compare the regularizing effects
of MINRES and MR-II, and shed light on some new features of them.
In Section \ref{SectionMR2}, we present our theoretical results on MR-II
with a detailed analysis. In Section \ref{SectionEx}, we numerically confirm
our theory that MINRES has only the partial regularization
for a general ill-posed problem and its hybrid variant is needed. Also, we
experimentally illustrate that MR-II has the full regularization for severely and
moderately ill-posed problems, which is stronger than our theory,
and it has the partial regularization for mildly ill-posed problems.
We also compare MR-II with LSQR, demonstrating that MR-II is as effective
as and at least twice as efficient as LSQR. We conclude the paper in
Section \ref{SectionCon}.

Throughout the paper, we denote by
$\mathcal{K}_{k}(C, w)=span\{w,Cw,\ldots,C^{k-1}w\}$
the $k$-dimensional Krylov subspace generated by the matrix
$\mathit{C}$ and the vector
$\mathit{w}$, by $\|\cdot\|$ and $\|\cdot\|_F$ the 2-norm of a matrix or vector
and the Frobenius norm of a matrix, respectively,
and by $I$ the identity matrix with order clear from the context.

\section{MINRES and MR-II}\label{SectionRev}

MINRES \cite{paige75} is based on the symmetric Lanczos process that
constructs an orthonormal basis of the Krylov subspace $\mathcal{K}_k(A,b)$.
Let $\bar{q}_1=b/\|b\|$. The $k$-step symmetric Lanczos process can be written
in the matrix form
\begin{equation*}
  A\bar{Q}_k=\bar{Q}_{k+1}\bar{T}_k,
\end{equation*}
where $\bar{Q}_{k+1}=(\bar{q}_1, \bar{q}_2, \ldots, \bar{q}_{k+1})$
has orthonormal columns which form $\mathcal{K}_k(A,b)$, and
$\bar{T}_k\in \mathbb{R}^{(k+1)\times k}$
is a tridiagonal matrix with its leading $k\times k$ submatrix symmetric.

At iteration $k$, MINRES solves $\|b-A\bar{x}^{(k)}\|=
\min_{x\in\mathcal{K}_k(A,b)}\|b-Ax\|$ for the iterate
$\bar{x}^{(k)}= \bar{Q}_k\bar{y}^{(k)}$ with
\begin{equation}\label{minres}
  \bar{y}^{(k)}=\arg\min\limits_{y\in \mathbb{R}^{k}}\|\|b\|e_1 - \bar{T}_ky\|,
\end{equation}
where $e_1$ is the first canonical vector of dimension $k+1$. For
our analysis purpose, it is important to write
\begin{equation}\label{minsol}
\bar{x}^{(k)}= \bar{Q}_k\bar{T}_k^{\dagger}\bar{Q}_{k+1}^Tb,
\end{equation}
which is the minimum 2-norm least squares solution of the perturbed
problem that replace $A$ in \eqref{eq1} by its rank $k$ approximation
$\bar{Q}_{k+1}\bar{T}_k\bar{Q}_k^T$.

MR-II \cite{hanke95} is
a variant of MINRES applied to $\mathcal{K}_k(A,Ab)$ which excludes the noisy $b$.
The method is based on the $k$-step symmetric Lanczos process
\begin{equation}\label{eqmform1}
  AQ_k=Q_{k+1}T_k,
\end{equation}
where $Q_{k+1}=(q_1, q_2, \ldots, q_{k+1})$ has orthonormal columns with
$q_1=Ab/\|Ab\|$,
$T_k\in \mathbb{R}^{(k+1)\times k}$ is a tridiagonal matrix
with the diagonals $\alpha_i$, the subdiagonals
$\beta_i>0,\ i=1,2,\ldots,k$ and the superdiagonals $\beta_i>0,\
i=1,2,\ldots,k-1$, and the first $k$ rows of $T_k$ is symmetric.
The columns of $Q_k$ form an orthonormal basis of $\mathcal{K}_k(A,Ab)$.
Mathematically, since the eigenvalues of $A$ are simple and $Ab$
has nonzero components in the directions
of all the eigenvectors $v_i$ of $A$, the Lanczos process can be run to $n$ steps
without breakdown, i.e., $\beta_k>0,\ k=1,2,\ldots,n-1$ and $\beta_n=0$.

At iteration $k$, MR-II solves $\|b-A x^{(k)}\|=
\min_{x\in\mathcal{K}_k(A,Ab)}\|b-Ax\|$ for the iterate
$x^{(k)}= Q_k y^{(k)}$ with
\begin{equation}\label{mrii}
  y^{(k)}=\arg\min\limits_{y\in \mathbb{R}^{k}}\|b - Q_{k+1}T_ky\|.
\end{equation}
Similar to \eqref{minsol}, we have the expression
\begin{equation}\label{mriisol}
x^{(k)}= Q_kT_k^{\dagger}Q_{k+1}^Tb,
\end{equation}
which is the minimum 2-norm least squares solution of the perturbed
problem that replace $A$ in \eqref{eq1} by its rank $k$ approximation
$Q_{k+1}T_kQ_k^T$.

The significance of \eqref{minsol} and \eqref{mriisol} lies that
the MINRES and MR-II iterates are the minimum 2-norm least squares
solutions of the perturbed
problems that replace $A$ in \eqref{eq1} by its rank $k$ approximations
$\bar{Q}_{k+1}\bar{T}_k\bar{Q}_k^T$ and $Q_{k+1}T_kQ_k^T$, respectively,
whose $k$ nonzero singular values are just those of $\bar{T}_k$ and $T_k$,
respectively. If the singular values of $\bar{T}_k$ or $T_k$ approximate
the $k$ large singular values of $A$ in natural order for $k=1,2,\ldots,k_0$,
then $\bar{Q}_{k+1}\bar{T}_k\bar{Q}_k^T$ and $Q_{k+1}T_kQ_k^T$ are near best
rank $k$ approximations to $A$ with accuracy similar to that of the best rank
$k$ approximation $A_k$. If this is the case,
MINRES and MR-II must have the full regularization, and
$\bar{x}^{(k_0)}$ and $x^{(k_0)}$ are best possible regularized solutions and
are as accurate as the best TSVD regularized solution $x_{k_0}$.

\section{The regularizing effects of MINRES}\label{SectionMin}

Similar to the CGLS and LSQR iterates \cite[p. 146]{hansen98}, we can
establish the following result on the MINRES iterates.

\begin{thm}\label{thm1}
For MINRES to solve \eqref{eq1} with the starting vector $\bar{q}_1=b/\|b\|$,
the $k$th iterate $\bar{x}^{(k)}$ has the form
\begin{equation}\label{eqfilter2}
  \bar{x}^{(k)}=\sum\limits_{i=1}^nf_i^{(k)}\frac{v_i^{T}b}{\lambda_i}v_i,
\end{equation}
where the filters $f_i^{(k)}=1-\prod\limits_{j=1}^k\frac{\theta_j^{(k)}-\lambda_i}
{\theta_j^{(k)}},\ i=1,2,\ldots,n$ with $|\lambda_1|> |\lambda_2|>\cdots>
|\lambda_n|>0$, and $\theta_j^{(k)},\ j=1,2,\ldots,k$, are the harmonic Ritz
values of $A$ with respect to $\mathcal{K}_k(A,b)$ and labeled as
$|\theta_1^{(k)}|> |\theta_2^{(k)}|>\cdots>|\theta_k^{(k)}|>0$.
\end{thm}

{\em Proof}.
From \cite{paige95}, the residual $\bar{r}^{(k)}=b-A\bar{x}^{(k)}$ of
the MINRES iterate $\bar{x}^{(k)}$ can be written as
\begin{equation}\label{res}
\bar{r}^{(k)}=\chi_k(A)b,
\end{equation}
where the residual polynomial $\chi_k(t)$ has the form
\begin{equation*}
  \chi_k(t)=\prod\limits_{j=1}^k
  \frac{\theta_j^{(k)}-t}{\theta_j^{(k)}}
\end{equation*}
with the $\theta_j^{(k)}$ the harmonic Ritz values of $A$ with respect to
$\mathcal{K}_k(A,b)$. From \eqref{res}, we get
\begin{equation*}
  \bar{x}^{(k)} = (I-\chi_k{(A)})A^{-1}b.
\end{equation*}
Substituting $A=V\Lambda V^T$ into the above gives
\begin{equation*}\label{}
  \bar{x}^{(k)}=\sum\limits_{i=1}^nf_i^{(k)}\frac{v_i^{T}b}{\lambda_i}v_i,
\end{equation*}
where
\begin{equation*}
  f_i^{(k)}=1-\prod\limits_{j=1}^k
  \frac{\theta_j^{(k)}-\lambda_i}{\theta_j^{(k)}},\ \ i=1,2,\ldots,n.
  \qed
\end{equation*}

Relation \eqref{eqfilter2} shows that the MINRES iterate $\bar{x}^{(k)}$
has a filtered SVD expansion. For a general symmetric $A$,
the harmonic Ritz values have an attractive feature: they usually
favor extreme eigenvalues of $A$, provided that a Krylov subspace contains
substantial information on all the eigenvectors $v_i$ \cite{paige95}.
In our current context, if at least a small
harmonic Ritz value in magnitude starts to appear for some $k\leq k_0$, i.e.,
$|\theta_{k}^{(k)}|\leq |\lambda_{k_0+1}|$, the corresponding filter factors
$f_i^{(k)}$, $i=k+1,\ldots,n$, are not small, meaning that $\bar{x}^{(k)}$
is already deteriorated. On the other hand, if no small
harmonic Ritz value in magnitude appears before $k\leq k_0$, the $\bar{x}^{(k)}$ are
expected to become better approximations to $x_{true}$ until $k=k_0$.
Unfortunately, since $\mathcal{K}_k(A,b)$
includes the noisy $b=\hat{b}+e$, which contains non-negligible components of
$v_i$ corresponding to small eigenvalues $\lambda_i$, it is generally possible
that a small harmonic Ritz value can appear for $k\leq k_0$.
This demonstrates that, in general, MINRES only has the partial regularization
and cannot obtain a best possible regularized solution.

\section{Regularization relationships between MINRES and MR-II}\label{SectionCom}

It was known a long time ago that MR-II has better regularizing effects
than MINRES, that is, MR-II obtains a more accurate regularized solution
than MINRES does \cite{hanke95}. Such phenomenon is simply due to
the fact that $\mathcal{K}_k(A,b)$ for MINRES includes the noisy
$b$ and $\mathcal{K}_k(A,Ab)$ for MR-II contains less
information on $v_i$ corresponding to small eigenvalues in magnitude
since the noise $e$ in the starting vector $Ab$ is filtered by
multiplication with $A$. Previously, we have given an analysis on the
regularizing effects of MINRES and shown that a hybrid MINRES is generally
needed for an ill-posed problem, independent of the degree of ill-posedness
of \eqref{eq1}. Next we shed more light on the regularization of
MINRES, compare it with MR-II, and reveal some new features of them.

To simplify our discussions, without loss of generality, we can well assume that
for a standard nonsingular linear system, the smaller residual, the more accurate
the approximate  solution is.  Given the residual minimization property of
MINRES and MR-II, one might be confused that, since
$\mathcal{K}_{k-1}(A,Ab)\subset \mathcal{K}_k(A,b)$,
the $k$th MINRES iterate $\bar{x}^{(k)}$ should be
at least as accurate as the $(k-1)$th MR-II iterate $x^{(k-1)}$. This is true
for solving the standard linear system where the right-hand
side is supposed to be {\em exact}, but it is nontrivial and depends for
solving an ill-posed problem, for which the $b$ is noisy and
we are concerned with
regularized approximations to the true solution $x_{true}$ other than the naive
solution $x_{naive}$. Our previous analysis has shown that a small harmonic
Ritz value $|\theta_k^{(k)}|\leq \sigma_{k_0+1}=|\lambda_{k_0+1}|$ generally
appears for MINRES before some iteration $k\leq k_0$, causing that MINRES
has only the partial regularization. On the other hand, however, note that
the regularized solutions $\bar{x}^{(k)}$ by MINRES converge to $x_{true}$ until
the semi-convergence of MINRES. As a result, because $\mathcal{K}_{k-1}(A,Ab)\subset
\mathcal{K}_k(A,b)$, the $k$th MINRES iterate $\bar{x}^{(k)}$ is more accurate than
the $(k-1)$th MR-II iterate $x^{(k-1)}$ {\em only until} the semi-convergence
of MINRES.

We can also explain the partial regularization of MINRES in terms of
the rank $k$ approximation $\bar{Q}_{k+1}\bar{T}_k\bar{Q}_k^T$ to $A$ as follows:
Since the $k$-dimensional dominant eigenspace of $A$ is identical to
its $k$-dimensional dominant left and right singular subspaces,
$\mathcal{K}_k(A,b)$ contains substantial information on all the $v_i$.
As a result, it is generally possible that the projected matrix $\bar{T}_k$ has
a singular value smaller than $\sigma_{k_0+1}$ for some $k\leq k_0$.
This means that $\bar{Q}_{k+1}\bar{T}_k\bar{Q}_k^T$ is a poor
rank $k$ approximation to $A$, causing, from \eqref{minres}, that
$\|\bar{x}^{(k)}\|=\|\bar{Q}_k\bar{y}^{(k)}\|=\|b\|\|\bar{T}_k^{\dagger}e_1\|$
is generally large, i.e., $\bar{x}^{(k)}$ is already deteriorated. Conversely,
if no singular value of $\bar{T}_k$ is smaller
than $\sigma_{k_0+1}$ and the semi-convergence of MINRES
does not yet occur, the MINRES iterate $\bar{x}^{(k)}$ should be at least as
accurate as the MR-II iterate $x^{(k-1)}$
because of $\mathcal{K}_{k-1}(A,Ab)\subset \mathcal{K}_k(A,b)$.

In summary, we need to use a hybrid MINRES with the TSVD method or
the standard-form Tikhonov regularization applied to the projected
problem in \eqref{minres} to expand the Krylov subspace until it contains
all the $k_0$ dominant spectral components and a best regularized
solution is found, in which the additional regularization
aims to remove the effects of small singular values of $\bar{T}_{k+1}$,
similar to the hybrid LSQR see \cite[Section 6.4]{hansen10}.

\section{Regularizing effects of MR-II}
\label{SectionMR2}

Before proceeding, we point out that, unlike \eqref{eqfilter2} for the MINRES iterates
$\bar{x}^{(k)}$, we have found that the MR-II iterates $x^{(k)}$ do not have filtered
SVD expansions of similar form. Even so, we can establish a number of other results
that help to better understand the regularization of MR-II. We first investigate
a fundamental problem: how well does the underlying subspace $\mathcal{K}_k(A,Ab)$
capture the $k$ dimensional dominant eigenspace of $A$? This problem is of basic
importance because it critically affects the accuracy of $Q_{k+1}T_kQ_k^T$ as a
rank $k$ approximation to $A$.

In terms of the definition of canonical angles
$\Theta(\mathcal{X},\mathcal{Y})$ between the two subspaces
$\mathcal{X}$ and $\mathcal{Y}$ of the same dimension \cite[p. 250]{stewart01},
we present the following result.

\begin{thm}\label{thm2}
Let $A=V\Omega\Sigma V^T=V\Lambda V^T$ be defined as \eqref{decom},
and assume that the singular values of $A$ are $\sigma_j=|\lambda_j|
=\mathcal{O}(e^{-\alpha j})$ with $\alpha>0$.
Let $\mathcal{V}_k=span\{V_k\}$ be the $k$-dimensional dominant spectral subspace
spanned by the columns of
$V_k=(v_1,v_2,\ldots,v_k)$, and $\mathcal{V}_k^s=\mathcal{K}_{k}(A, Ab)$. Then
\begin{equation}\label{sindelta}
 \|\sin\Theta(\mathcal{V}_k,\mathcal{V}_k^s)\|=
 \frac{\|\Delta_k\|}{\sqrt{1+\|\Delta_k\|^2}}
\end{equation}
with the $(n-k)\times k$ matrix $\Delta_k$ to be defined by \eqref{deltak} and
\begin{equation}\label{eqres1}
 \|\Delta_k\|_F\leq
  \frac{|\lambda_{k+1}|}{|\lambda_{k}|}\frac{\max_{j=k+1}^n
  |v_{j}^Tb|}{\min_{j=1}^k|v_{j}^Tb|}
  \sqrt{k(n-k)}\left(1+\mathcal{O}(e^{-\alpha})\right), \ k=1,2,\ldots,n-1.
\end{equation}
\end{thm}

{\em Proof}.
Note that $\mathcal{K}_{k}(\Lambda,\Lambda V^Tb)$ is spanned by the columns
of the $n\times k$ matrix $DB_k$
with
\begin{equation*}\label{}
  D=\mathrm{diag}\left(\lambda_i v_i^Tb\right),\ \ \
  B_k=\left(\begin{array}{cccc} 1 &
  \lambda_1&\ldots & \lambda_1^{k-1}\\
1 &\lambda_2 &\ldots &\lambda_2^{k-1} \\
\vdots & \vdots&&\vdots\\
1 &\lambda_n &\ldots &\lambda_n^{k-1}
\end{array}\right).
\end{equation*}
Partition $D$ and $B_k$ as follows:
\begin{equation*}\label{}
  D=\left(\begin{array}{cc} D_1 & 0 \\ 0 & D_2 \end{array}\right),\ \ \
  B_k=\left(\begin{array}{c} B_{k1} \\ B_{k2} \end{array}\right),
\end{equation*}
where $D_1, B_{k1}\in\mathbb{R}^{k\times k}$. Since $B_{k1}$ is
a Vandermonde matrix with $\lambda_j$ distinct for $1\leq j\leq k$, it is
nonsingular. Noting $\mathcal{K}_k(A,Ab)=V\mathcal{K}_{k}
(\Lambda,\Lambda V^Tb)$, we have
\begin{equation*}\label{}
  \mathcal{K}_{k}(A, Ab)=span\{VDB_k\}=span
  \left\{V\left(\begin{array}{c} D_1B_{k1} \\ D_2B_{k2} \end{array}\right)\right\}
  =span\left\{V\left(\begin{array}{c} I\\ \Delta_k \end{array}\right)\right\}
\end{equation*}
with
\begin{equation}\label{deltak}
\Delta_k=D_2B_{k2}B_{k1}^{-1}D_1^{-1}.
\end{equation}
Define $Z_k=V\left(\begin{array}{c} I \\ \Delta_k \end{array}\right)$.
Then $Z_k^TZ_k=I+\Delta_k^T\Delta_k$, and the columns of
$Z_k(Z_k^TZ_k)^{-\frac{1}{2}}$ form an orthonormal basis of $\mathcal{V}_k^s$.

Write $V=(V_k, V_k^{\perp})$. By definition, we obtain
\begin{align}
   \|\sin\Theta(\mathcal{V}_k,\mathcal{V}_k^s)\|  &=
   \left\|(V_k^{\perp})^T Z_k(Z_k^TZ_k)^{-\frac{1}{2}}\right\|\notag\\
   &=\left\|(V_k^{\perp})^T V
   \left(\begin{array}{c} I \\ \Delta_k \end{array}\right)
   (I+\Delta_k^T\Delta_k)^{-\frac{1}{2}}\right\| \notag\\
   &=\|\Delta_k(I+\Delta_k^T\Delta_k)^{-1/2}\|\notag\\
   &=\frac{\|\Delta_k\|}{\sqrt{1+\|\Delta_k\|^2}},\notag
\end{align}
which proves \eqref{sindelta}.

We next estimate $\|\Delta_k\|$ and establish upper bound for the right-hand
side of \eqref{sindelta}.  We have
\begin{align}
\|\Delta_k\|&\leq\|\Delta\|_F=\left\|D_2B_{k2}B_{k1}^{-1}D_1^{-1}\right\|_F
   \leq \|D_2\|\left\|B_{k2}B_{k1}^{-1}\right\|_F\left\|D_1^{-1}\right\|\notag\\
 &=\frac{|\lambda_{k+1}|}{|\lambda_{k}|}\frac{\max_{j=k+1}^n|v_{j}^Tb|}
 {\min_{j=1}^k|v_{j}^Tb|}\left\|B_{k2}B_{k1}^{-1}\right\|_F.
   \label{fnorm}
\end{align}

We now estimate $\left\|B_{k2}B_{k1}^{-1}\right\|_F$.
It is easily justified that the $i$th column of $B_{k1}^{-1}$ consists of
the coefficients of the Lagrange polynomial
\begin{equation*}\label{}
  L_i^{(k)}(\lambda)=\prod\limits_{j=1,j\neq i}^k
  \frac{\lambda-\lambda_j}{\lambda_i-\lambda_j}
\end{equation*}
that interpolates the elements of the $i$th canonical basis vector
$e_i^{(k)}\in \mathbb{R}^{k}$ at the abscissas $\lambda_1,
\ldots, \lambda_k$. Consequently, the $i$th column of $B_{k2}B_{k1}^{-1}$ is
\begin{equation*}\label{}
  B_{k2}B_{k1}^{-1}e_i^{(k)}=\left(L_i^{(k)}(\lambda_{k+1}),
  \ldots,L_i^{(k)}(\lambda_{n})\right)^T,
\end{equation*}
from which we obtain
\begin{equation}\label{tk2tk1}
  B_{k2}B_{k1}^{-1}=\left(\begin{array}{cccc} L_1^{(k)}(\lambda_{k+1})&
  L_2^{(k)}(\lambda_{k+1})&\ldots & L_k^{(k)}(\lambda_{k+1})\\
L_1^{(k)}(\lambda_{k+2})&L_2^{(k)}(\lambda_{k+2}) &\ldots &
L_k^{(k)}(\lambda_{k+2}) \\
\vdots & \vdots&&\vdots\\
L_1^{(k)}(\lambda_{n})&L_2^{(k)}(\lambda_{n}) &\ldots &L_k^{(k)}(\lambda_{n})
\end{array}\right).
\end{equation}
For a fixed $\lambda$ satisfying $|\lambda|\leq |\lambda_{k+1}|$,
let $i_0=\arg\displaystyle\max_{i=1,2,\ldots,k}|L_i^{(k)}(\lambda)|$.
Then we have
\begin{align}\label{li0k}
   |L_{i_0}^{(k)}(\lambda)|&=\prod\limits_{j=1,j\neq i_0}^k
  \left|\frac{\lambda-\lambda_j}{\lambda_{i_0}-\lambda_j}\right|
  \leq \prod\limits_{j=1,j\neq i_0}^k
  \left|\frac{|\lambda_j-\lambda|}{|\lambda_j|-|\lambda_{i_0}|}\right|
  \leq \prod\limits_{j=1,j\neq i_0}^k
  \left|\frac{|\lambda_j|+|\lambda_{k+1}|}{|\lambda_j|-|\lambda_{i_0}|}\right|.
\end{align}
Therefore, for $i=1,2\ldots,k$ and $|\lambda|\leq |\lambda_{k+1}|$,
making use of Taylor series expansions, we get
\begin{align}\label{}
|L_i^{(k)}(\lambda)| & \leq\prod\limits_{j=1,j\neq i_0}^k
\left|\frac{|\lambda_j|+|\lambda_{k+1}|}{|\lambda_j|-|\lambda_{i_0}|}\right|
  =\prod\limits_{j=1}^{i_0-1}\frac{|\lambda_j|+
   |\lambda_{k+1}|}{|\lambda_j|-|\lambda_{i_0}|}
   \cdot\prod\limits_{j=i_0+1}^{k}\frac{|\lambda_j|+
   |\lambda_{k+1}|}{|\lambda_{i_0}|-|\lambda_{j}|}\notag\\
& =\prod\limits_{j=1}^{i_0-1}\frac{1+\mathcal{O}\left(e^{-(k-j+1)\alpha}\right)}
{1-\mathcal{O}\left(e^{-(i_0-j)\alpha}\right)}\cdot
\prod\limits_{j=i_0+1}^{k}\frac{\mathcal{O}\left(e^{-(k-j+1)\alpha}\right)+1}
{\mathcal{O}\left(e^{(j-i_0)\alpha}\right)-1} \notag\\
&=\frac{\prod\limits_{j=1}^{k}\left(1+\mathcal{O}\left(e^{-(k-j+1)\alpha}\right)
\right)}
  {1+\mathcal{O}\left(e^{-(k-i_0+1)\alpha}\right)}
\prod\limits_{j=1}^{i_0-1}\frac{1}
{1-\mathcal{O}(e^{-(i_0-j)\alpha})}
\prod\limits_{j=i_0+1}^{k}\frac{1}
{\mathcal{O}(e^{(j-i_0)\alpha})-1}\notag\\
&=\frac{\prod\limits_{j=1}^{k}\left(1+\mathcal{O}\left(e^{-(k-j+1)\alpha}\right)\right)}
  {(1+\mathcal{O}\left(e^{-(k-i_0+1)\alpha}\right))}
\prod\limits_{j=1}^{i_0-1}\frac{1}{1-\mathcal{O}(e^{-(i_0-j)\alpha})} \notag\\
  & \quad \quad\cdot
\prod\limits_{j=i_0+1}^{k}\frac{1}{1-\mathcal{O}(e^{-(j-i_0)\alpha})}
\frac{1}{\prod\limits_{j=i_0+1}^{k}\mathcal{O}\left(e^{(j-i_0)\alpha}\right)}\notag\\
&= \frac{\left(1+\sum\limits_{j=1}^{k+1}\mathcal{O}\left(e^{-(k-j+1)\alpha}\right)\right)}
  {(1+\mathcal{O}\left(e^{-(k-i_0+1)\alpha}\right))}
\frac{\left(1+\sum\limits_{j=1}^{i_0} \mathcal{O}(e^{-j\alpha})\right)
\left(1+\sum\limits_{j=1}^{k-i_0+1} \mathcal{O}(e^{-j\alpha})\right)}
{\prod\limits_{j=i_0+1}^{k}\mathcal{O}(e^{(j-i_0)\alpha})} \label{lik}
\end{align}
by absorbing those higher order terms into ${\cal O}(\cdot)$.
Note that in the above numerator we have
  $$
  1+\sum\limits_{j=1}^{k+1}\mathcal{O}(e^{-(k-j+1)\alpha})
  =1+\mathcal{O}\left(\sum\limits_{j=1}^{k+1}e^{-(k-j+1)\alpha}\right)
  =1+ \mathcal{O}\left(\frac{e^{-\alpha}}
    {1-e^{-\alpha}}(1-e^{-(k+1)\alpha})\right),
  $$
  $$
  1+\sum\limits_{j=1}^{i_0} \mathcal{O}(e^{-j\alpha})
    =1+ \mathcal{O}\left(\sum\limits_{j=1}^{i_0}e^{-j\alpha}\right)
    =1+ \mathcal{O}\left(\frac{e^{-\alpha}}
    {1-e^{-\alpha}}(1-e^{-i_0\alpha})\right),
  $$
  and
  $$
    1+\sum\limits_{j=1}^{k-i_0+1} \mathcal{O}(e^{-j\alpha})
    =1+ \mathcal{O}\left(\sum\limits_{j=1}^{k-i_0+1}e^{-j\alpha}\right)
    =1+ \mathcal{O}\left(\frac{e^{-\alpha}}{1-e^{-\alpha}}
    (1-e^{-(k-i_0+1)\alpha})\right).
  $$
It is easy to check that for any $1\leq i_0\leq k$ the product of
the above three terms is no more than
  $$
  1+ \mathcal{O}\left(\frac{3e^{-\alpha}}{1-e^{-\alpha}}\right)
  +\mathcal{O}\left(\left(\frac{e^{-\alpha}}{1-e^{-\alpha}}\right)^2\right)
  =  1+\mathcal{O}(e^{-\alpha}).
  $$
By definition, the factor
$\prod\limits_{j=i_0+1}^{k}\mathcal{O}(e^{(j-i_0)\alpha})=
\prod\limits_{j=i_0+1}^{k}\frac{|\lambda_{i_0}|}{|\lambda_j|}$ in
the denominator of \eqref{lik}, which is exactly one when $i_0=k$, and it
is bigger than one when $i_0<k$; the other factor
$1+\mathcal{O}\left(e^{-(k-i_0+1)\alpha}\right)$ is between
$1+\mathcal{O}\left(e^{-k\alpha}\right)$ and
$1+\mathcal{O}\left(e^{-\alpha}\right)$.
Therefore, for any $k$ and $|\lambda|\leq |\lambda_{k+1}|$, we have
\begin{align}
|L_k^{(k)}(\lambda)|&=1+\mathcal{O}(e^{-\alpha}), \label{lkk}\\
  |L_{i_0}^{(k)}(\lambda)|&=\max_{i=1,2,\ldots,k} |L_i^{(k)}(\lambda)|=
  1+\mathcal{O}(e^{-\alpha}).\label{li0k}
\end{align}
From this estimate and \eqref{tk2tk1} it follows that
\begin{equation}\label{b1b2}
\left\|B_{k2}B_{k1}^{-1}\right\|_F\leq
\sqrt{k (n-k)}\left(1+\mathcal{O}(e^{-\alpha})\right).
\end{equation}
As a result, for $k=1,2,\ldots,n-1$, from \eqref{fnorm} we have
\begin{align*}
 \|\Delta\|_F
  &\leq\frac{|\lambda_{k+1}|}{|\lambda_{k}|}
   \frac{\max_{j=k+1}^n|v_{j}^Tb|}{\min_{j=1}^k|v_{j}^Tb|}
   \sqrt{k(n-k)}\left(1+\mathcal{O}(e^{-\alpha})\right). \qed
\end{align*}

{\bf Remark 5.1}\ \
Trivially, we have
$$
 \|\sin\Theta(\mathcal{V}_k,\mathcal{V}_k^s)\|\leq 1.
$$
But in our context it is impossible to have
$\|\sin\Theta(\mathcal{V}_k,\mathcal{V}_k^s)\|=1$ since
$\Delta_k$ is not a zero matrix. We have seen from the proof that the factor
$\frac{|\lambda_{k+1}|}{|\lambda_{k}|}
\frac{\max_{j=k+1}^n|v_{j}^Tb|}{\min_{j=1}^k|v_{j}^Tb|}$
in it is intrinsic and unavoidable in \eqref{eqres1}.
But the factor $\sqrt{n(n-k)}$ in \eqref{eqres1} is an overestimate
and can certainly be reduced. The reason is that \eqref{b1b2} is an overestimate
since $|L_i^{(k)}(\lambda_j)|$ for $i$ not near to $k$ is considerably
smaller than $|L_{i_0}^{(k)}(\lambda_j)|$,
$j=k+1,\ldots,n$ but we replace all them by the maximum
$|L_{i_0}^{(k)}(\lambda_j)|=1+\mathcal{O}(e^{-\alpha})$.
In fact, our derivation before \eqref{lkk} and \eqref{li0k}
when replacing $i_0$ by $i$
clearly illustrates that the smaller $i$ is, the smaller
$|L_i^{(k)}(\lambda_j)|$ than $|L_k^{(k)}(\lambda_j)|$, $j=k+1,\ldots,n$.

Recall the discrete Picard condition \eqref{picard}. Then the coefficients
\begin{equation}\label{ck}
  c_k=\frac{\max_{j=k+1}^n|v_{j}^Tb|}{\min_{j=1}^k|v_{j}^Tb|}=
  \frac{\max_{j=k+1}^n(|v_{j}^T\hat{b}+v_{j}^Te|)}
  {\min_{j=1}^k(|v_{j}^T\hat{b}+v_{j}^Te|)}\approx
\frac{|\lambda_{k+1}|^{1+\beta}+|v_{k+1}^Te|}{|\lambda_{k}|^{1+\beta}+|v_{k}^Te|}.
\end{equation}
We see that, the larger $\beta$ is, the smaller
$c_k\approx \frac{|\lambda_{k+1}|^{1+\beta}} {|\lambda_k|^{1+\beta}}$,
which is a constant for $k\leq k_0$, and thus the better
$\mathcal{V}_k^s$ captures $\mathcal{V}_k$. For $k>k_0$, since all the
$|v_k^T b|\approx |v_{k}^Te|$ are roughly the same, we have $c_k\approx 1$,
meaning that $\mathcal{V}_k^s$ may not capture $\mathcal{V}_k$ so well after
iteration $k_0$.

{\bf Remark 5.2}\ \ The theorem can be extended to
moderately ill-posed problems with the singular values
$\sigma_j=\mathcal{O}(j^{-\alpha})$, $\alpha>1$ considerably and $k$ not big,
where the factor $1+\mathcal{O}(e^{-\alpha})$ in \eqref{eqres1}
is replaced by a bigger $\mathcal{O}(1)$. Let us look into why it is so.
Recall that, by definition, $|L_{i_0}^{(k)}(\lambda)|\geq |L_k^{(k)}(\lambda)|$
for $|\lambda|\leq |\lambda_{k+1}|$.
Using a similar proof to that of Theorem~\ref{thm2} and the first order
Taylor expansion, we can roughly estimate
$|L_k^{(k)}(\lambda)|$ as follows:
\begin{align*}\label{}
|L_{i_0}^{(k)}(\lambda)|\approx |L_k^{(k)}(\lambda)|&\leq \prod\limits_{j=1}^{k-1}
\left|\frac{|\lambda_j|+|\lambda_{k+1}|}{|\lambda_j|-|\lambda_k|}\right| \\
&=\prod\limits_{j=1}^{k-1}\frac{1+\mathcal{O}\left((\frac{j}{k+1})^{\alpha}\right)}
{1-\mathcal{O}\left((\frac{j}{k})^{\alpha}\right)}\\
&\approx\sum\limits_{j=1}^{k-1} \left(1+\mathcal{O}\left(\left(\frac{j}
{k+1}\right)^{\alpha}\right)\right)
\cdot\sum\limits_{j=1}^{k-1} \left(1+\mathcal{O}\left(\left(\frac{j}
{k}\right)^{\alpha}\right)\right)
=\mathcal{O}(1).
\end{align*}
This estimate is not as accurate as that for severely ill-posed problems.
More important is that it depends on $k$ and
increases slowly as $k$ increases. The above estimate can be improved
when $A$ is symmetric definite:
\begin{align*}\label{}
|L_{i_0}^{(k)}(\lambda)|\approx |L_k^{(k)}(\lambda)|&= \prod\limits_{j=1}^{k-1}
\frac{|\lambda_j-\lambda_{k+1}|}{|\lambda_j-\lambda_k|}
\leq\prod\limits_{j=1}^{k-1}
\frac{|\lambda_j|}{|\lambda_j-\lambda_k|} \\
&=\prod\limits_{j=1}^{k-1}\frac{1}
{1-\mathcal{O}\left((\frac{j}{k})^{\alpha}\right)}\\
&\approx\sum\limits_{j=1}^{k-1} \left(1+\mathcal{O}\left(\left(\frac{j}
{k}\right)^{\alpha}\right)\right)
=\mathcal{O}(1),
\end{align*}
smaller than the previous one. The two estimates mean that $\mathcal{V}_k^s$ may
capture $\mathcal{V}_k$ better for $A$ symmetric definite than for $A$ symmetric
indefinite where there are both positive and negative ones among
the first $k+1$ large eigenvalues.

{\bf Remark 5.3}\ \ A combination of \eqref{sindelta} and \eqref{eqres1}
and the above analysis indicate that $\mathcal{V}_k^s$  captures
$\mathcal{V}_k$ better for severely ill-posed problems than for moderately
ill-posed problems. There are two reasons for this.
The first is that the factors $|\lambda_{k+1}/\lambda_k|$ are
basically fixed constants for severely ill-posed problems as $k$
increases, and they are smaller than the counterparts for moderately ill-posed
problems unless the degree $\alpha$ of its ill-posedness is far bigger than one
and $k$ small. The second is that the factor $1+\mathcal{O}(e^{-\alpha})$
is smaller for severely ill-posed problems than the factor
$\mathcal{O}(1)$ for moderately ill-posed problems.

{\bf Remark 5.4} \ \
The situation is fundamentally different for mildly ill-posed problems:
Firstly, we always have $|L_{i_0}^{(k)}(\lambda)|>1$ substantially for $|\lambda|
\leq |\lambda_{k+1}|$, $\alpha\leq 1$ and any $k$. Secondly,
$c_k$ defined by \eqref{ck} is closer to one than that
for moderately ill-posed problems for $k=1,2,\ldots,k_0$.
Thirdly, for the same noise level $\|e\|$ and $\beta$, from the discrete
Picard condition \eqref{picard} and the definition of $k_0$ we see that
$k_0$ is bigger for a mildly ill-posed problem than that for a moderately
ill-posed problem. All of them show that $\mathcal{V}_k^s$
captures $\mathcal{V}_k$ {\em considerably better} for severely and moderately
ill-posed problems than for mildly ill-posed problems.
In other words, our results imply that $\mathcal{V}_k^s$ contains substantial
information on the other $n-k$ eigenvectors for mildly ill-posed problems, causing
that a small harmonic Ritz value generally appears for some $k\leq k_0$,
especially when $k_0$ is not small. Equivalently, the projected matrix
$T_k$ generated by MR-II generally has a small singular value for some
$k\leq k_0$, such that the solution $x^{(k)}$ is deteriorated, as deduced from
\eqref{mriisol}. As a result, we are certain that MR-II has better regularizing
effects for severely and moderately ill-posed problems than for mildly ill-posed
problems. Most importantly, by this property, since MR-II has {\em at most} the
full regularization for severely and moderately ill-posed problems, we
deduce and are thus sure that MR-II generally has only the partial regularization
for mildly ill-posed problems.

We mention that, in comparison with the results, i.e., Theorem 2.1,
in \cite{huang14} on LSQR, we find that $\mathcal{K}_k(A,Ab)$ is as
comparably effective as $\mathcal{K}_k(A^TA,A^Tb)$, on which LSQR works,
for capturing the $k$-dimensional dominant eigenspace.

Let us get more insight into the regularization of MR-II.
Recall \eqref{mriisol}, where
$$
x^{(k)}=(Q_{k+1}T_kQ_k^T)^{\dagger}b
=Q_kT_k^{\dagger}Q_{k+1}^Tb.
$$
Define
\begin{equation}\label{defgamma}
\gamma_k = \left\|A-Q_{k+1}T_kQ_k^T\right\|,
\end{equation}
which measures the quality or accuracy of the rank $k$ approximation $Q_{k+1}T_kQ_k^T$
to $A$. This quantity is central and fundamental to understand the regularizing
effects of MR-II
and measures how the iterates $x^{(k)}$ by MR-II behave like the TSVD regularized
solution $x_k=A_k^{\dagger}b$. Particularly, note that the best rank $k_0$ approximation
$A_{k_0}$ satisfies $\|A-A_{k_0}\|=\sigma_{k_0+1}$. Then if
$\gamma_{k_0}\approx \sigma_{k_0+1}$ for  $\sigma_{k_0+1}$ reasonably small,
$Q_{k_0+1}T_{k_0}Q_{k_0}^T$ is a near best rank $k_0$ approximation to $A$
with approximate accuracy $\sigma_{k_0+1}$ and has no small nonzero singular value.
In this case, the regularized solution $x^{(k_0)}$ is close to the best TSVD regularized
solution $x_{k_0}$, and MR-II has the full regularization. Otherwise,
if $\gamma_{k_0}>\sigma_{k_0+1}$ considerably, then $Q_{k_0+1}T_{k_0}Q_{k_0}^T$
deviates from the best rank $k$ approximation $A_{k_0}$ considerably and
$x^{(k_0)}$ is not close to $x_{k_0}$, meaning that MR-II has only the partial
regularization.

Based on Theorem~\ref{thm2} and Remark 5.2, we can derive the following estimates
for $\gamma_k$.

\begin{thm}\label{thm3}
Assume that \eqref{eq1} is severely or moderately ill posed. Then
\begin{equation}\label{gamma}
|\lambda_{k+1}| \leq\gamma_k\leq |\lambda_{k+1}| +
|\lambda_1|\|\sin\Theta(\mathcal{V}_k,\mathcal{V}_k^s)\|.
\end{equation}
\end{thm}

{\em Proof}.
Note that $Q_{k+1}T_kQ_k^T$ is of rank $k$. The lower bound in
\eqref{gamma} is trivial since the best $k$ approximation $A_k$
to $A$ satisfies $\|A-A_k\|=\sigma_{k+1}=|\lambda_{k+1}|$.
We next prove the upper bound. From \eqref{eqmform1}, we obtain
\begin{align}
\left\|A-Q_{k+1}T_kQ_k^T\right\|&= \left\|A-AQ_kQ_k^T\right\|
=\left\|A(I-Q_kQ_k^T)\right\|.\label{newform}
\end{align}
From Theorem \ref{thm2}, it is known that $\mathcal{V}_k^s
=\mathcal{K}_{k}(A,Ab)=span\{Q_k\}$.
Let $V_k=(v_1, v_2, \ldots, v_k)$ and
$\Lambda_k=\mathrm{diag}(\lambda_1,\lambda_2,\ldots,\lambda_k)$.
Then by the definition of $\|\sin\Theta(\mathcal{V}_k,\mathcal{V}_k^s)\|$
we obtain
\begin{align*}
 \left\|A-AQ_kQ_k^T\right\| &= \left\|(A-V_k\Lambda_kV_k^T+V_k\Lambda_kV_k^T)
     (I-Q_kQ_k^T)\right\| \notag\\
 &\leq\left\|(A-V_k\Lambda_kV_k^T)(I-Q_kQ_k^T)\right\|+
     \left\|V_k\Lambda_kV_k^T(I-Q_kQ_k^T)\right\| \notag\\
 &\leq |\lambda_{k+1}| + \|\Lambda_k\|\left\|V_k^T(I-Q_kQ_k^T)\right\| \notag\\
 &= |\lambda_{k+1}| + |\lambda_1|\|\sin\Theta(\mathcal{V}_k,\mathcal{V}_k^s)\|.
  \qed
\end{align*}

Our later numerical experiments
will indicate that $\gamma_k\approx \sigma_{k+1}=|\lambda_{k+1}|$
for severely and moderately ill-posed problems, illustrating that
$Q_{k+1}T_kQ_k^T$ is a near best rank $k$ approximation to $A$ with
the approximate accuracy $\sigma_{k+1}$. Particularly, since
$\gamma_{k_0}\approx \sigma_{k_0+1}$, the MR-II iterate
$x^{(k_0)}=Q_{k_0}T_{k_0}^{\dagger}Q_{k_0+1}^Tb$ is close
to the TSVD solution $x_{k_0}$ provided that $\sigma_{k_0+1}$.
Furthermore, we will find that the error
$\|x^{(k_0)}-x_{true}\|$ of MR-II iterate $x^{(k_0)}$
is as small as the error $\|x_{k_0}-x_{true}\|$
of the best TSVD solution $x_{k_0}$. This indicates that MR-II has the full
regularization. Experimentally, for severely and moderately ill-posed problems,
the observations $\gamma_k\approx \sigma_{k+1}$ appear to be general
and thus should have strong theoretical supports. Our upper bound
in \eqref{gamma} appears to be a considerable overestimate.

Recall that $\alpha_i$ and $\beta_i,\,i=1,2,\ldots,k$ denote the
diagonals and subdiagonals of $T_k$ defined by \eqref{eqmform1},
respectively.
We next establish some interesting and intimate relationships
between them and $\gamma_k$, showing how fast $\alpha_k$ and $\beta_k$
decay.

\begin{thm}\label{thm4}
For $k=1,2,\ldots,n-2$ we have
\begin{eqnarray}
  \beta_{k+1}&\leq& \gamma_k,\label{alphak1}\\
    |\alpha_{k+2}|&\leq &\gamma_k. \label{betak}
\end{eqnarray}
\end{thm}

{\em Proof}.
Since the Lanczos process can be run to completion, we have
\begin{equation*}
Q_n^T  A Q_n=\hat{T}_n,
\end{equation*}
where $Q_n \in\mathbb{R}^{n\times n}$ is orthogonal, and
\begin{equation}\label{hattn}
   \hat{T}_n = \left(\begin{array}{ccccc}
   \alpha_1 & \beta_1& & &\\
   \beta_1 & \alpha_2 & \beta_2 & & \\
   &\ddots & \ddots &\ddots & \\
   & & &\ddots & \beta_{n-1} \\
   & & &\beta_{n-1} & \alpha_{n}
   \end{array}\right)
\end{equation}
is symmetric tridiagonal. Thus, from \eqref{eqmform1} we have
\begin{align*}
  \gamma_k &= \left\|A-Q_{k+1}T_kQ_k^T\right\|
           = \left\|Q_n^T(A-Q_{k+1}T_kQ_k^T)Q_n\right\| \\
           &=\left\|\hat{T}_n- \left(\begin{array}{c} I\\
           \mathbf{0} \end{array}\right)T_k
           \left(\begin{array}{cc} I & \mathbf{0}
           \end{array}\right)\right\|
          =\|G_k\|,
\end{align*}
where
\begin{equation*}
  G_k = \left(\begin{array}{ccccc}
   \beta_{k+1} & \alpha_{k+2}&\beta_{k+2} & &\\
   &\beta_{k+2} & \alpha_{k+3}  &\beta_{k+4} & \\
   & &\ddots & \ddots  & \\
   & & &\ddots & \beta_{n-1} \\
   & & &\beta_{n-1} & \alpha_{n}
   \end{array}\right) \in \mathbb{R}^{(n-k-1)\times (n-k)},
\end{equation*}
from which and $\beta_i>0,\ i=1,2,\ldots,n-1$ it follows that
$$
\beta_{k+1}=\|G_ke_1\|\leq\|G_k\|=\gamma_k
$$
and
\begin{equation*}
|\alpha_{k+2}|\leq\sqrt{\alpha_{k+2}^2+\beta_{k+2}^2}=\|G_ke_2\|\leq\|G_k\|=
 \gamma_k
\end{equation*}
for $k=1,2,\ldots,n-2$.
Therefore, \eqref{alphak1} and \eqref{betak} hold. \qed

This theorem indicates that $|\alpha_{k+2}|$ and $\beta_{k+1}$ decay at least
as fast as $\gamma_k$. Moreover, based on the experimental observations that
$\gamma_k\approx\sigma_{k+1}$ for severely and moderately ill-posed
problems, the theorem illustrates that $|\alpha_{k+2}|$ and
$\beta_{k+1}$ decay as fast as $\sigma_{k+1}$,
$k=1,2,\ldots,n-2$, for these two kinds of problems.

\section{Numerical experiments}\label{SectionEx}

In this section, we report numerical experiments to illustrate the
regularizing effects of MINRES and MR-II and make a number of comparisons.
We justify our theory: (i) MINRES has only the partial regularization,
independent of the degree of ill-posedness, and a hybrid MINRES
is generally needed;
(ii) the $k$th MINRES iterate $\bar{x}^{(k)}$ is always more accurate than
the $(k-1)$th MR-II iterate $x^{(k-1)}$ until the semi-convergence of MINRES;
(iii) MR-II has only the partial regularization for mildly ill-posed
problems, and a hybrid MR-II is needed. In the meantime, experimentally, we
demonstrate a remarkable and attractive feature, stronger than our theory predicts:
MR-II has the full regularization
for severely and moderately ill-posed problems and its iterates at semi-convergence
is as accurate as the best TSVD solutions for these two kinds of problems.
We will use the function ${\sf lcurve}$ in \cite{hansen07} to depict the L-curves.
In order to simulate exact arithmetic, the Lanczos process with reorthogonalization
is used in MINRES and MR-II.

Table~\ref{table} lists test problems and their degree of ill-posedness,
all of which are symmetric and arise from the
discretization of the first kind Fredholm integral equations; see Hansen's
regularization toolbox \cite{hansen07} for details. For each problem except
the 2D image deblurring problem `blur', we use the corresponding
code in \cite{hansen07} to
generate a $1024\times 1024$ matrix $A$, the true solution $x_{true}$
and noise-free right-hand $\hat{b}$. In order to simulate the noisy data,
we generate the white noise vector $e$ whose entries are normally
distributed with mean zero, such that the relative noise level
$\varepsilon=\frac{\|e\|}{\|\hat{b}\|}=10^{-2},10^{-3}, 10^{-4}$, respectively.
To simulate exact arithmetic, the full reorthogonalization is used
during the Lanczos process. We remind that, as far as ill-posed problem \eqref{eq1}
is concerned,
our primary goal consists in justifying the regularizing effects
of iterative solvers, which are {\em unaffected by sizes}
of ill-posed problems and only depends on the degree of
ill-posedness. Therefore, for this purpose, as extensively done in the
literature (see, e.g., \cite{hansen98,hansen10} and the references therein),
it is enough to test not very large problems. Indeed, for
$n$ large, say, 1,0000, we have observed completely the same
behavior as that for $n$ not large, e.g., $n=1024$ used in this paper except for
the problem `blur' with $n=65,536$. A reason for using $n$ not large is because
such choice makes it practical to fully
justify the regularization effects of LSQR by comparing it with
the TSVD method, which suits only for
small and/or medium sized problems for computational efficiency.
All the computations are carried out using Matlab 7.8 with the machine
precision $\epsilon_{\rm mach}=2.22\times10^{-16}$ under the Microsoft
Windows 7 64-bit system.
\begin{table}[h]
    \caption{The description of test problems.}
    \centering
   \begin{tabular*}{12.1cm}{@{\extracolsep{\fill}}lll}
     \hline\smallskip
     Problem & Description & Ill-posedness \\
     \hline
     shaw & One-dimensional image restoration model & severe\\
     foxgood & Severely ill-posed test problem & severe\\
     gravity & One-dimensional gravity surveying problem& severe\\
     phillips & phillips' "famous" test problem & moderate\\
     deriv2 & Computation of second derivative & mild\\
     blur &   2D Image deblurring test problem & mild/moderate \\
     \hline\smallskip
   \end{tabular*}
   \label{table}
  \end{table}

\subsection{A comparison of the regularizing effects of MR-II and MINRES}
\label{seccom}

We now compare MINRES and MR-II and justify our theory: (i)
the MR-II iterate is always more accurate
than the MINRES iterate at their respective semi-convergence, meaning that
MINRES cannot obtain best possible regularized solutions
and has only the partial regularization, independent of the degree of
ill-posedness; (ii) the MINRES iterates $\bar{x}^{(k)}$ are always
more accurate that the MR-II iterates $x^{(k-1)}$ until the semi-convergence of
MINRES.

In this subsection, we only report the results for the noise level
$\varepsilon=10^{-3}$. Results for the other two $\varepsilon$ are analogous
and thus omitted.

\begin{figure}[t]

\begin{minipage}{0.48\linewidth}
  \centerline{\includegraphics[width=7.0cm,height=5cm]{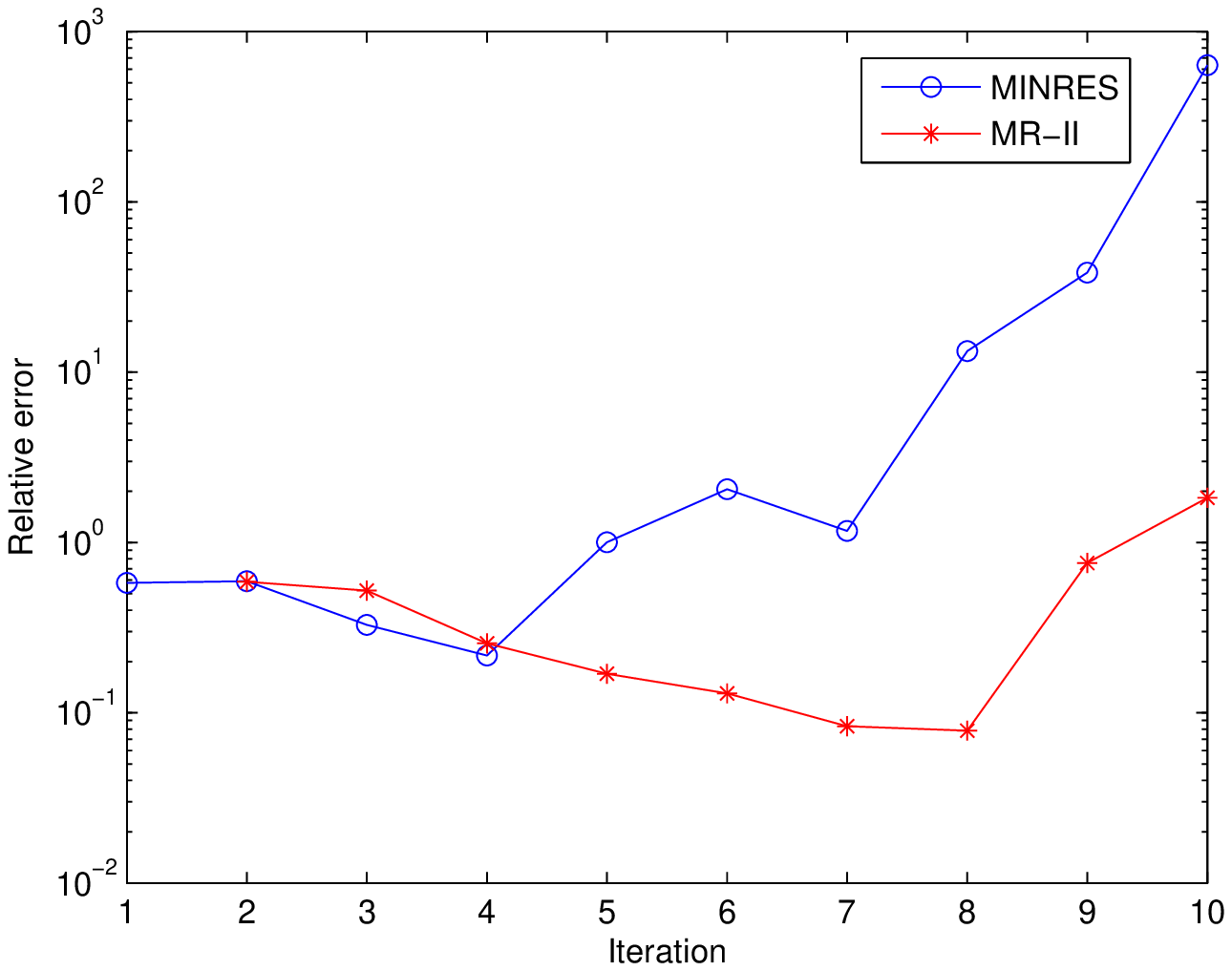}}
  \centerline{(a)}
\end{minipage}
\hfill
\begin{minipage}{0.48\linewidth}
  \centerline{\includegraphics[width=7.0cm,height=5cm]{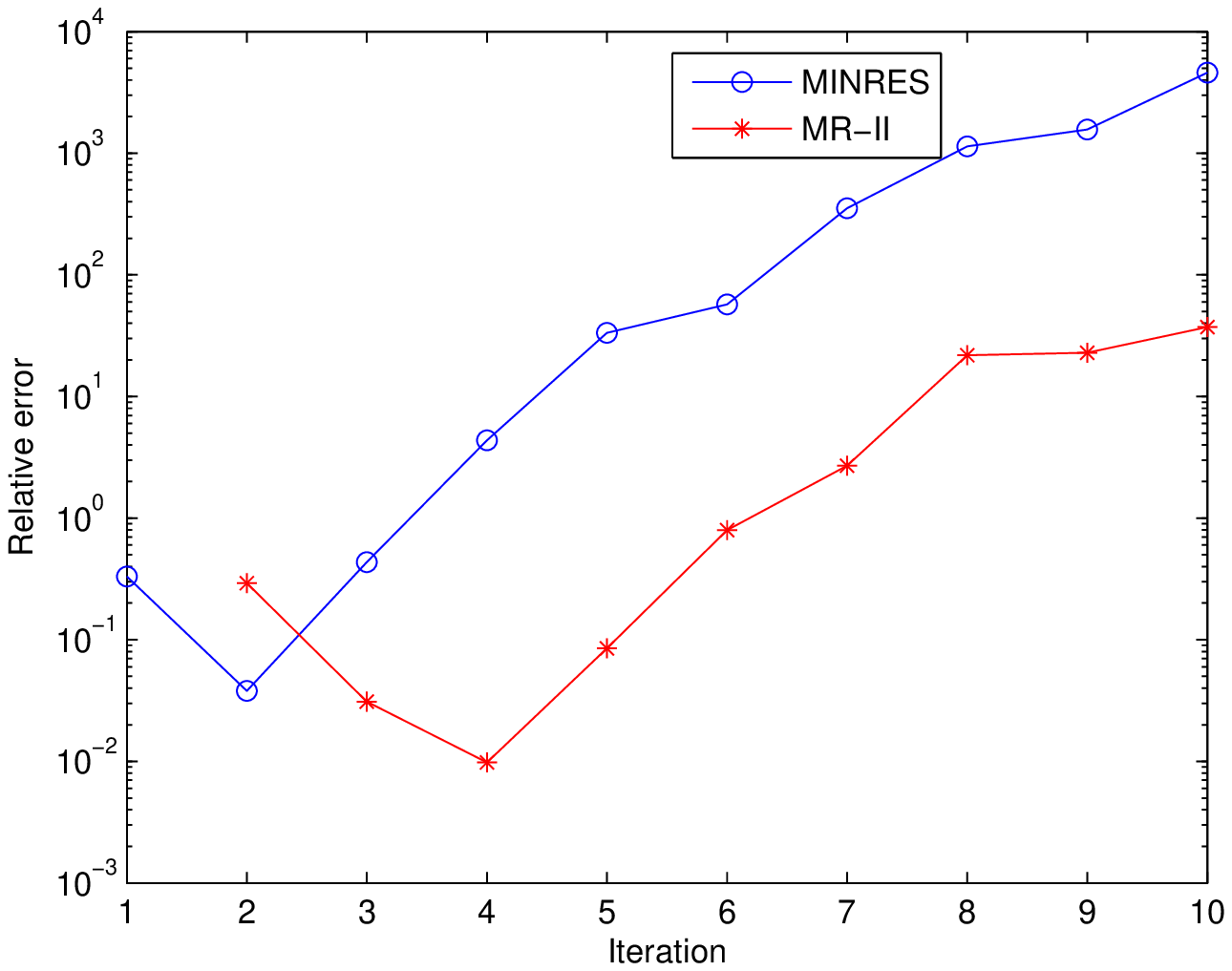}}
  \centerline{(b)}
\end{minipage}

\vfill
\begin{minipage}{0.48\linewidth}
  \centerline{\includegraphics[width=7.0cm,height=5cm]{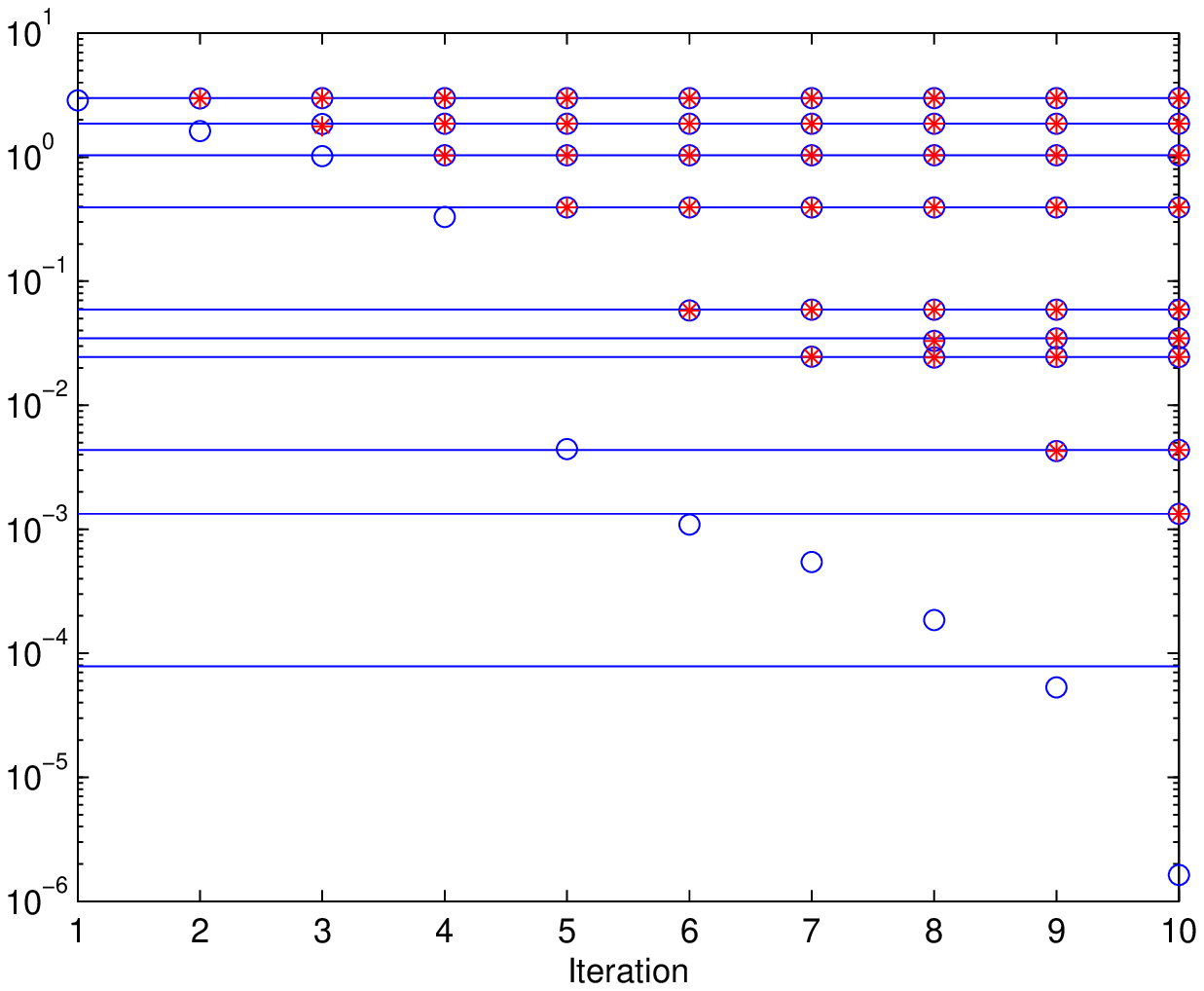}}
  \centerline{(c)}
\end{minipage}
\hfill
\begin{minipage}{0.48\linewidth}
  \centerline{\includegraphics[width=7.0cm,height=5cm]{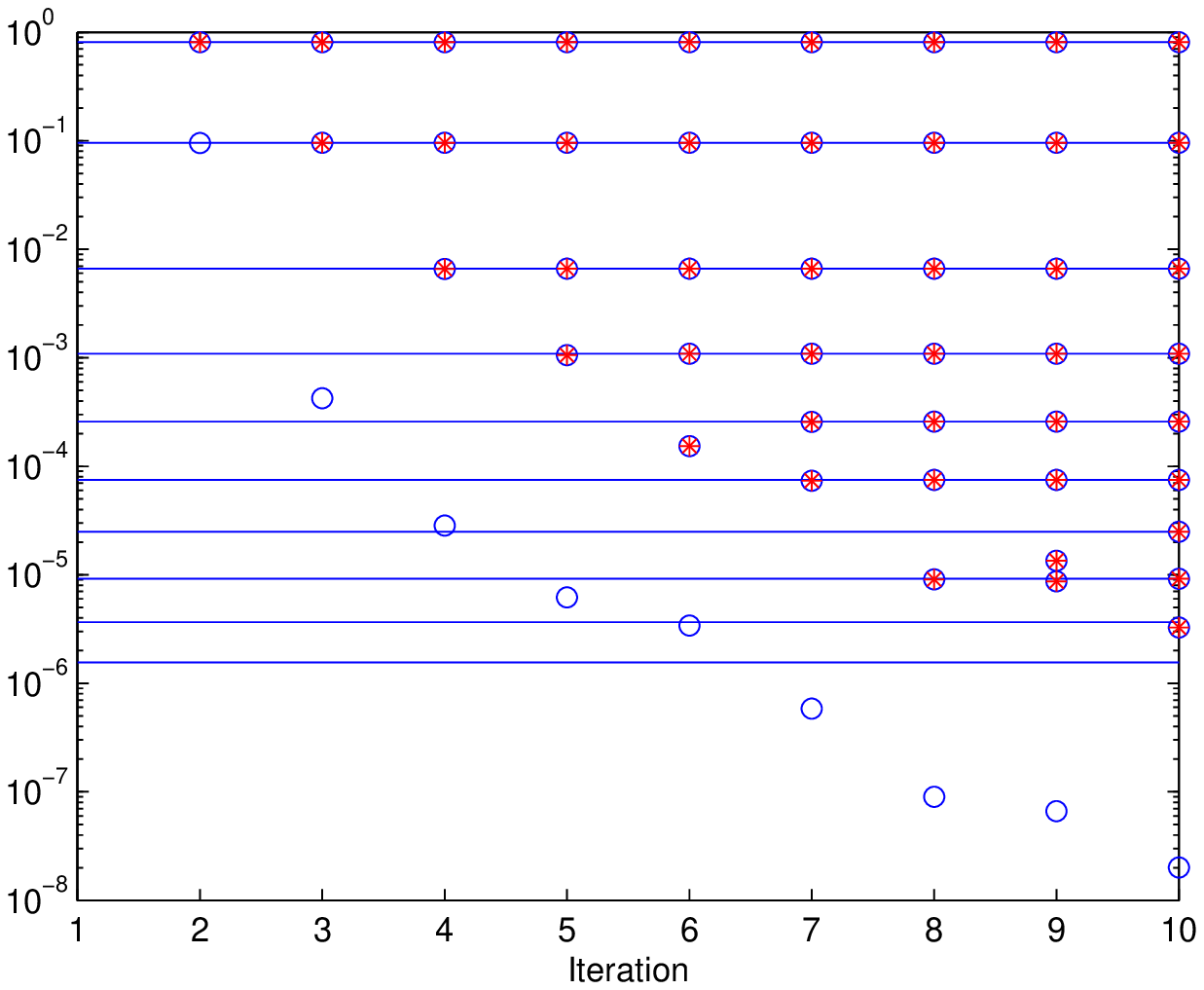}}
  \centerline{(d)}
\end{minipage}
\caption{(a)-(b): The relative errors $\|x^{(k)}-x_{true}\|/\|x_{true}\|$
by MINRES and MR-II;
(c)-(d): Plots of the singular values (circles for MINRES, stars for
MR-II) of
the projected matrices and the ones (solid lines) of $A$ for
shaw (left) and foxgood (right).}
\label{fig1}
\end{figure}

\begin{figure}[t]

\begin{minipage}{0.48\linewidth}
  \centerline{\includegraphics[width=7.0cm,height=5cm]{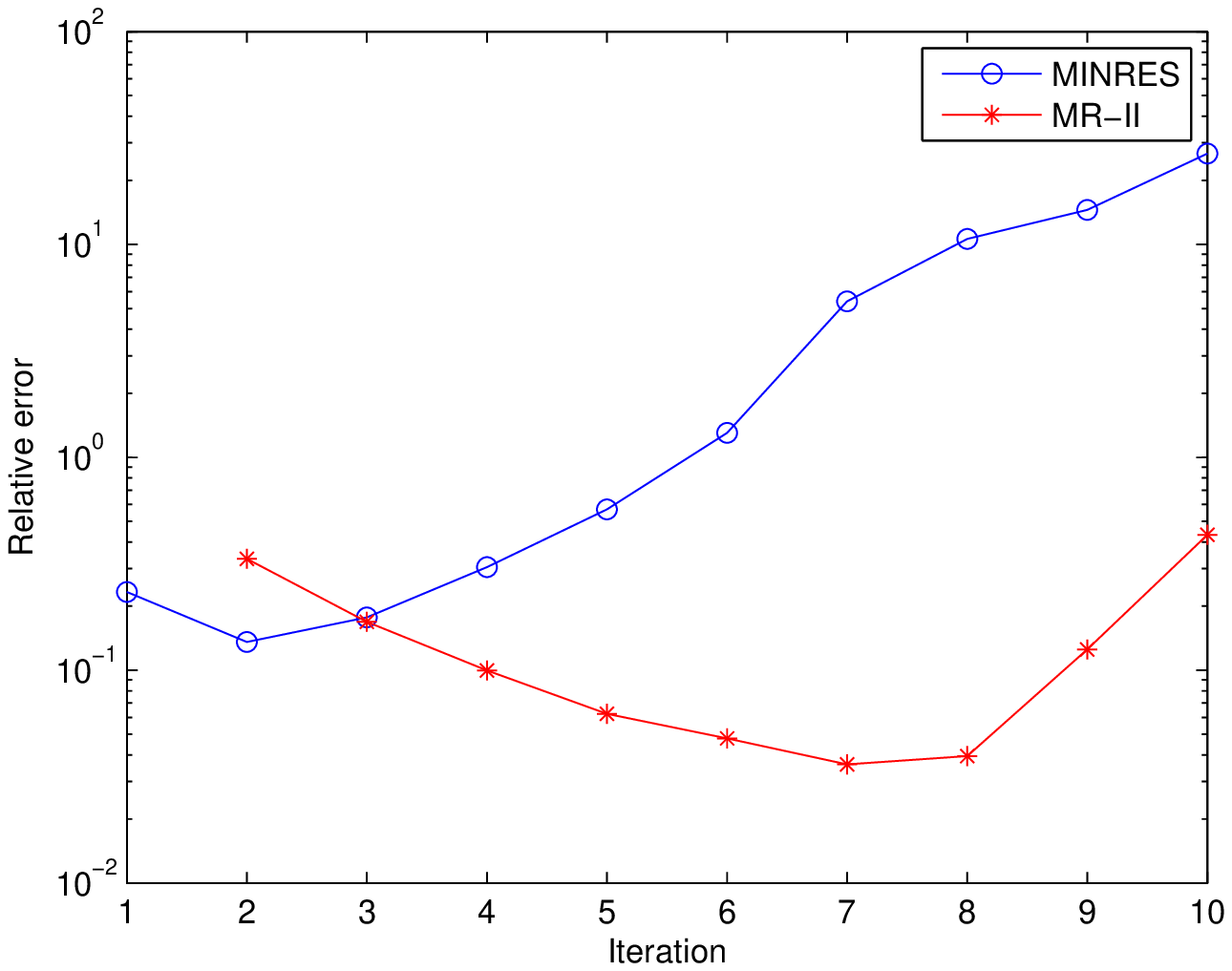}}
  \centerline{(a)}
\end{minipage}
\hfill
\begin{minipage}{0.48\linewidth}
  \centerline{\includegraphics[width=7.0cm,height=5cm]{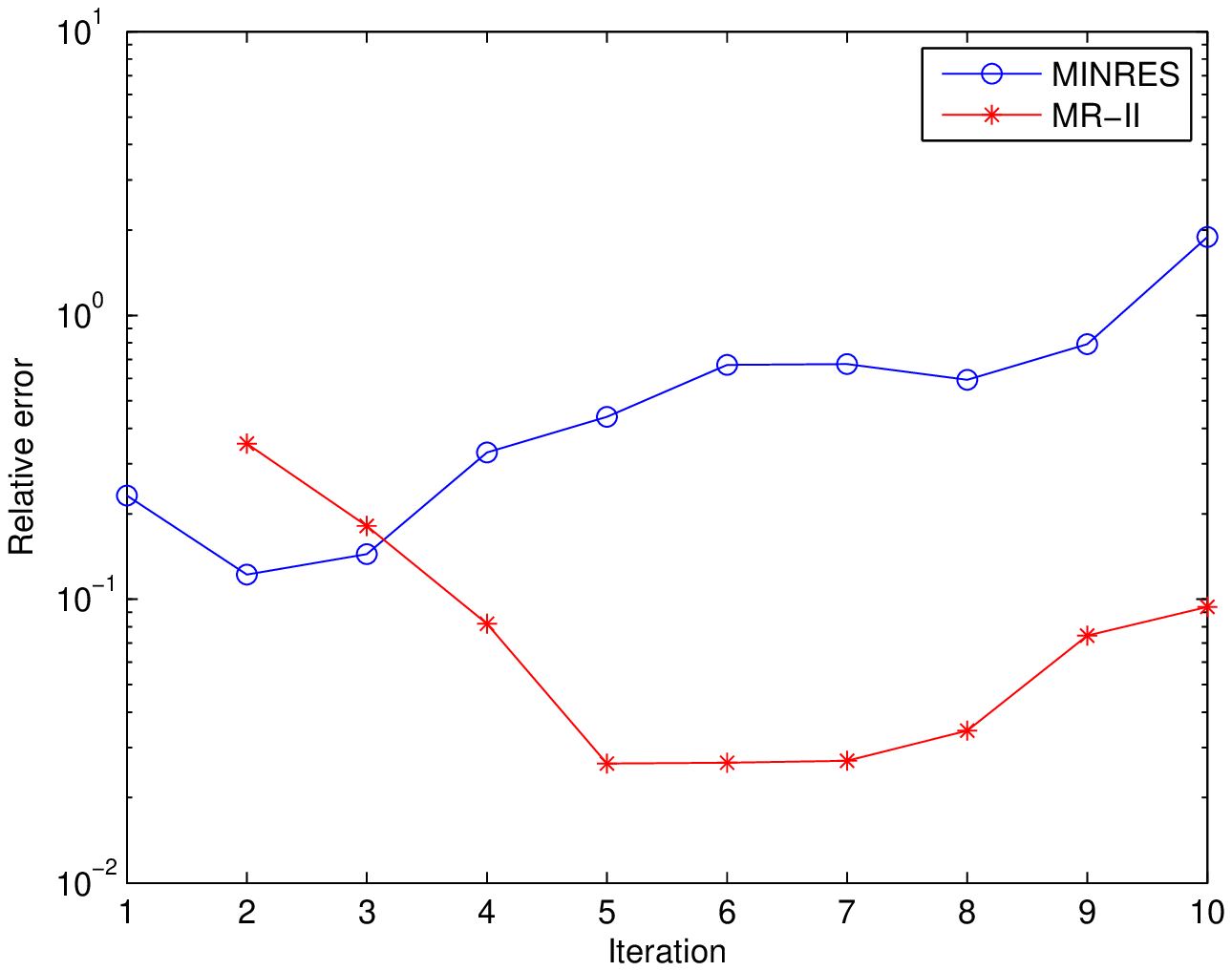}}
  \centerline{(b)}
\end{minipage}
\vfill
\begin{minipage}{0.48\linewidth}
  \centerline{\includegraphics[width=7.0cm,height=5cm]{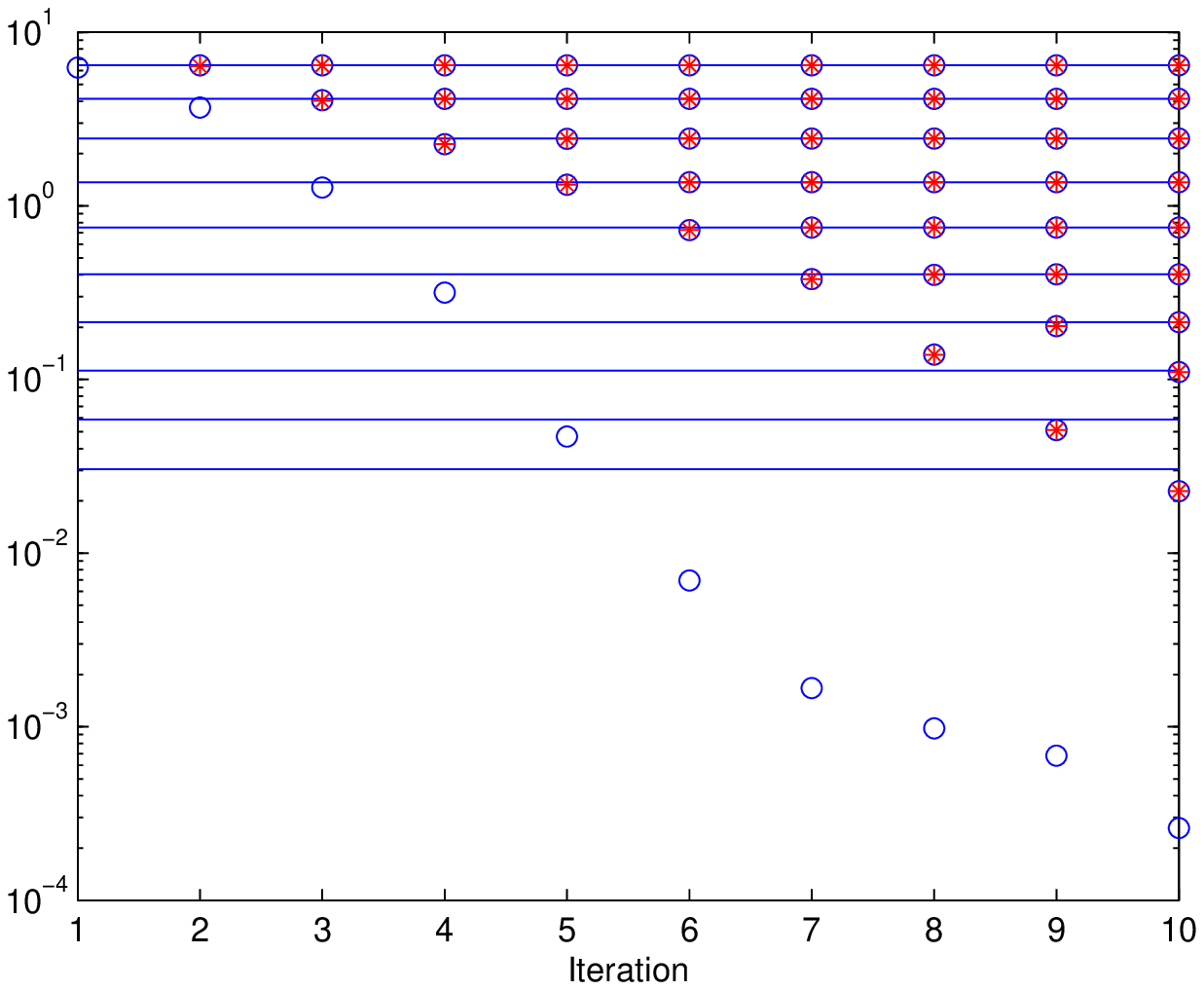}}
  \centerline{(c)}
\end{minipage}
\hfill
\begin{minipage}{0.48\linewidth}
  \centerline{\includegraphics[width=7.0cm,height=5cm]{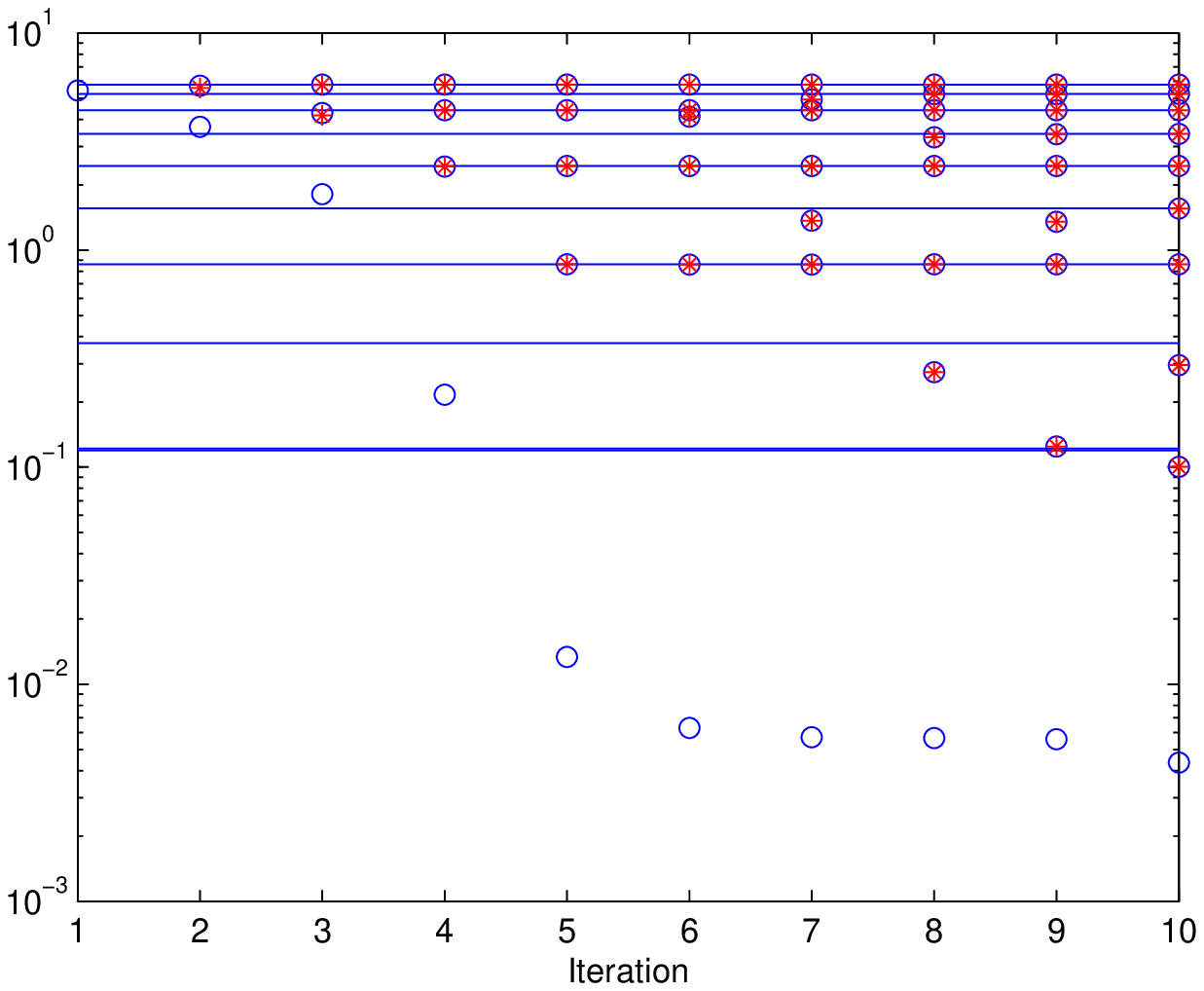}}
  \centerline{(d)}
\end{minipage}
\caption{ (a)-(b): The relative errors $\|x^{(k)}-x_{true}\|/\|x_{true}\|$
by MINRES and MR-II;
(c)-(d): Plots of the singular values (circles for MINRES, stars for
MR-II) of the projected matrices
and the ones (solid lines) of $A$ for gravity (left) and phillips (right). }
\label{fig2}
\end{figure}

Figures~\ref{fig1} and \ref{fig2} display numerous curves for severely and
moderately ill-posed problems. Clearly, all the
MR-II ierates are always more accurate than the MINRES iterates at their
respective semi-convergence. This indicates that MINRES has
only the partial regularization. As elaborated previously,
this is because that a small singular value of the projected matrix $\bar{T}_k$
appears before a regularized solution becomes best,
causing that its error does not reach the same error level as that obtained by
MR-II. For instance, we see from
Figure~\ref{fig1} (a) and (c) that all the singular values of $T_k$ in MR-II
are excellent approximations to the $k$ large singular values of $A$ in natural
order for $k\leq 9$. We see that the semi-convergence of MR-II occurs at
iteration $k=7$. By the comments in the end of Section~\ref{SectionRev} and
the explanations after \eqref{defgamma},
this clearly justifies the full regularization of MR-II, and the best possible
regularized solution by MR-II includes {\em seven} dominant spectral or SVD components.
On the other hand, it is clearly seen from Figure~\ref{fig1} (c) that
the smallest singular value of $\bar{T}_5$ in MINRES is smaller than $\sigma_8=\sigma_{k_0+1}$,
making the relative error starts to increase dramatically at iteration 5 and MINRES
have only the partial regularization.

Similar phenomena are observed for foxgood, and MR-II has the full regularization
with $k_0=3$, as indicated by
Figure~\ref{fig1} (b) and (d), where the smallest singular value
of $\bar{T}_3$ lies between $\sigma_4$ and $\sigma_5$ and the best iterate
$\bar{x}^{(k)}$ by MINRES at semi-convergence is considerably less accurate
than the best iterate $x^{(k)}$ by MR-II at semi-convergence,
meaning that MINRES has only the partial regularization.
We have analogous findings for gravity and phillips,
as shown by Figure~\ref{fig2} (b) and (d), which again demonstrate that MR-II has
the full regularization but MINRES has only the partial regularization.

As for the mildly ill-posed problem deriv2, we also see from Figure~\ref{fig3} (a)
that the relative error obtained by
MR-II clearly reaches the lower minimum level than that by MINRES, indicating
that MR-II has better regularizing effects than MINRES.

The above experiments have illustrated
that MR-II always obtains more accurate regularized solutions
than MINRES does for the test severely, moderately and mildly problems.
This justifies our theory that MINRES only has the
partial regularization, independent of the degree of ill-posedness. Therefore,
one must use a hybrid MINRES with some regularization method applied to the
projected problems in order to remove the effects of small singular values of
$\bar{T}_k$ and improve the accuracy of regularized solutions until a best
regularized solution is found.

It is clear from Figures~\ref{fig1}--\ref{fig2} and Figure~\ref{fig3} (a) that, for
each test problem, the MINRES iterates $\bar{x}^{(k)}$ are more accurate than
the corresponding MR-II iterates $x^{(k-1)}$ until the
semi-convergence of MINRES. Afterwards,
the regularized solutions $\bar{x}^{(k)}$ are deteriorated more and more
seriously. This confirms our theory in Section~\ref{SectionCom},
i.e., assertion (ii) in the beginning of this subsection.

\subsection{The regularizing effects of MR-II, MINRES and their
hybrid variants}
\label{Sectionfullregu}

We first test MR-II, MINRES and their hybrid variants for the
mildly ill-posed problem deriv2, and justify our theory that MR-II
has only the partial regularization and one must use its hybrid
variant to compute a best possible regularized solution.

\begin{figure}

\begin{minipage}{0.48\linewidth}
  \centerline{\includegraphics[width=7.0cm,height=5cm]{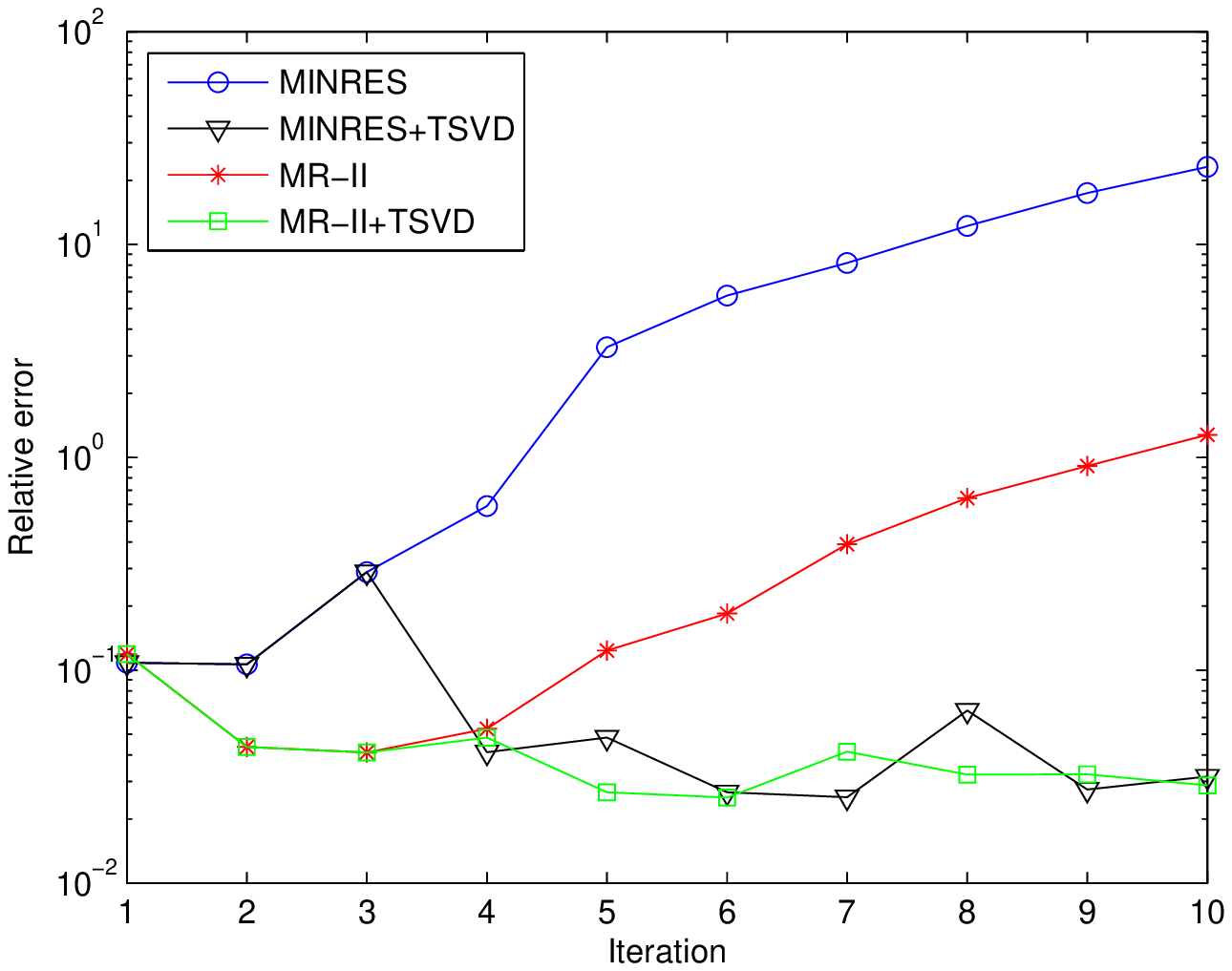}}
  \centerline{(a)}
\end{minipage}
\hfill
\begin{minipage}{0.48\linewidth}
  \centerline{\includegraphics[width=7.0cm,height=5cm]{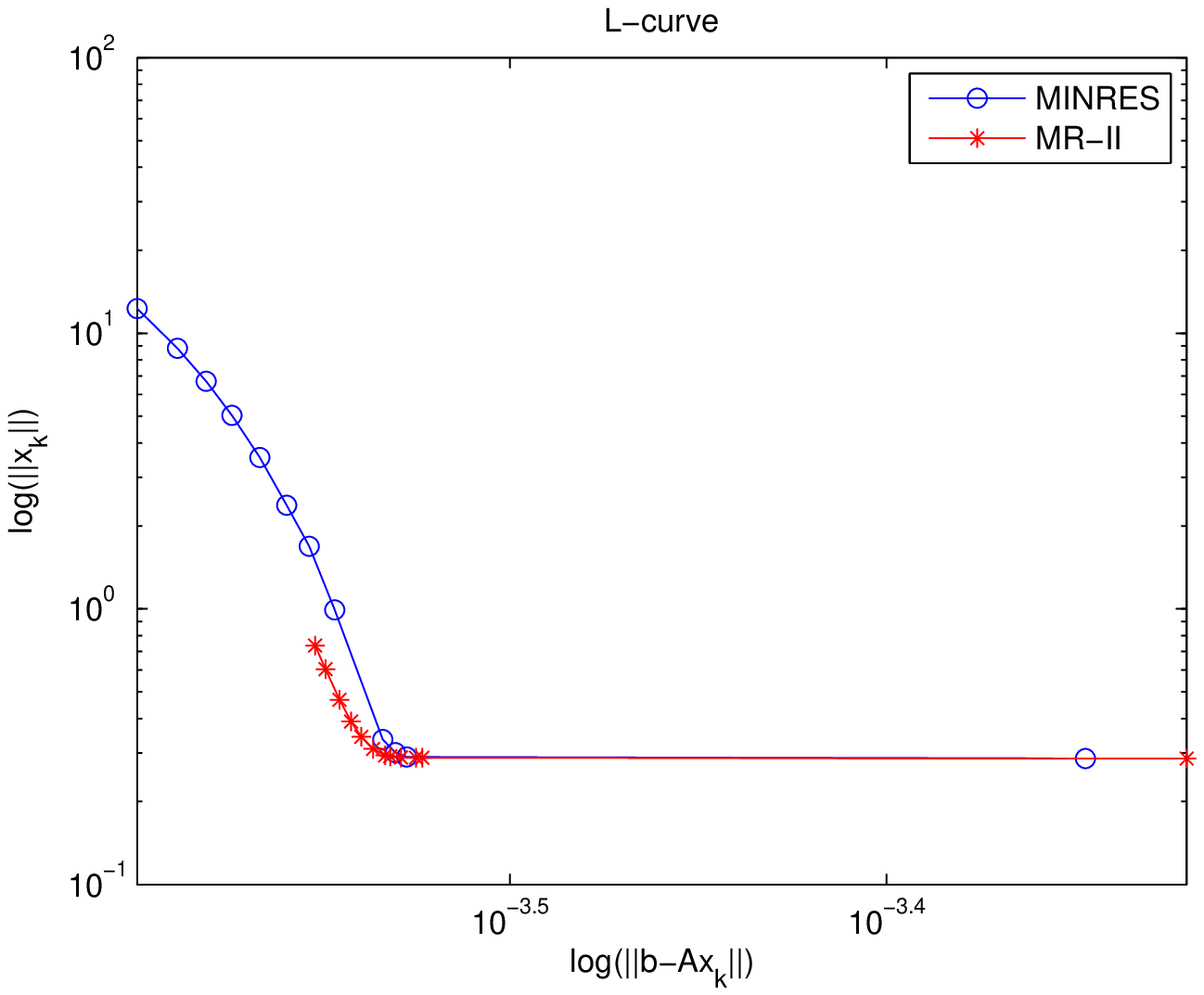}}
  \centerline{(b)}
\end{minipage}
\caption{(a): The relative errors $\|x^{(k)}-x_{true}\|/\|x_{true}\|$
by the pure MINRES and MR-II as well as the hybrid MINRES and MR-II; (b): The L-curves
of MINRES and MR-II for deriv2.}
\label{fig3}
\end{figure}

For deriv2, Figure~\ref{fig3} (a) shows that the relative errors of regularized
solutions obtained by the hybrid MINRES and MR-II with the TSVD regularization
method applied to the projected problems reach
a considerably smaller minimum level than those by MINRES and MR-II themselves.
For this problem, before MINRES or MR-II
captures all the dominant spectral components needed,
a small singular value of $\bar{T}_k$ or $T_k$ appears and starts to deteriorate
the regularized solutions. In contrast, their hybrid variants
expand Krylov subspaces until enough dominant spectral components are captured
and the TSVD regularization method effectively dampens the SVD components
corresponding to small singular values of the projected matrices
$\bar{T}_k$ by MINRES and $T_k$ by MR-II. For example, we see from
Figure~\ref{fig3} (a) that the semi-convergence of
MR-II occurs at iteration $k=3$,
which is also observed by the corner of the L-curve depicted by
Figure~\ref{fig3} (b). However, as shown by Figure~\ref{fig3} (a),
such regularization of MR-II is not enough, and
the hybrid MR-II uses a larger
six dimensional Krylov subspace $\mathcal{K}_6(A,Ab)$ to improve
the solutions and get a best possible regularized solution, whose residual
norm is smaller than that obtained by the pure
MR-II. After $k=6$,
the regularized solutions almost stabilize with the minimum error
as $k$ increases. We observe similar phenomena for MINRES and its
hybrid variant, where we find that the relative error by the hybrid
MINRES reaches the same minimum level as that by the hybrid
MR-II.

Next we test MR-II, MINRES and their hybrid variants for severely and
moderately ill-posed problems. We attempt
to get more insight into the regularizing effects of MR-II.
As a matter of fact, we have already justified the full regularization
of MR-II for the four test problems in Section~\ref{seccom}. In what follows,
we will give more details and justifications on
the full regularization of MR-II. We show that (i) the relative error obtained by
the hybrid MINRES reaches the same minimum level as that by the hybrid MR-II;
(ii) MR-II has the full regularization effects:
at semi-convergence, the regularized solution by the pure MR-II is as accurate
as that by the
hybrid MR-II with the TSVD regularization used within projected problems;
(iii) MR-II generates near best rank $k$ approximations
$Q_{k+1}T_kQ_k^T$ to $A$, i.e., the relation $\gamma_k\approx\sigma_{k+1}=|\lambda_{k+1}|$
holds with different noise levels. Keep in mind \eqref{tsvdsolution} and \eqref{mriisol}.
This means that $Q_{k+1}T_kQ_k^T$ generated by MR-II plays the same role as $A_k$, the
best rank $k$ approximation to $A$, so that MR-II has the full regularization.

\begin{figure}[t]

\begin{minipage}{0.48\linewidth}
  \centerline{\includegraphics[width=7.0cm,height=5cm]{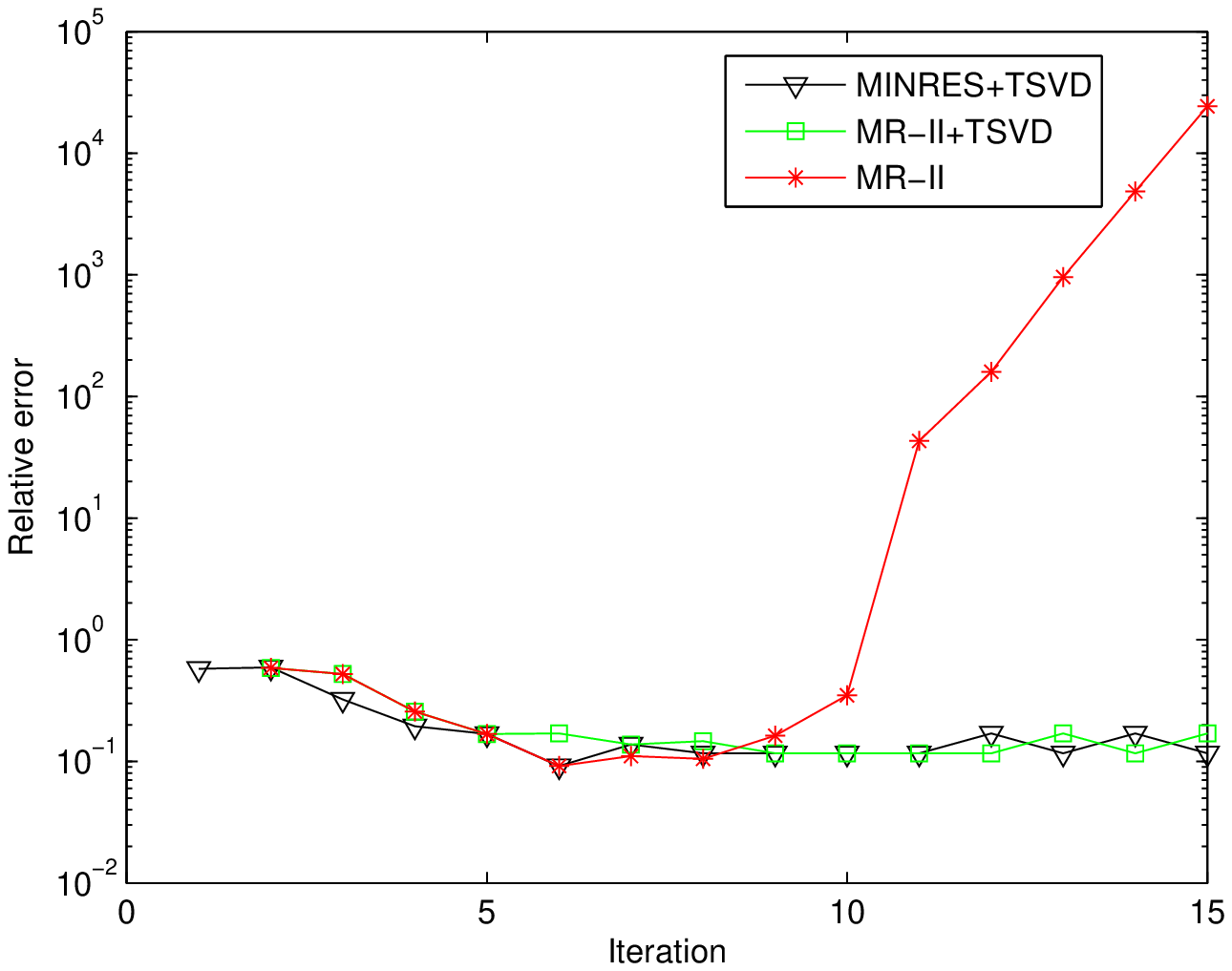}}
  \centerline{(a)}
\end{minipage}
\hfill
\begin{minipage}{0.48\linewidth}
  \centerline{\includegraphics[width=7.0cm,height=5cm]{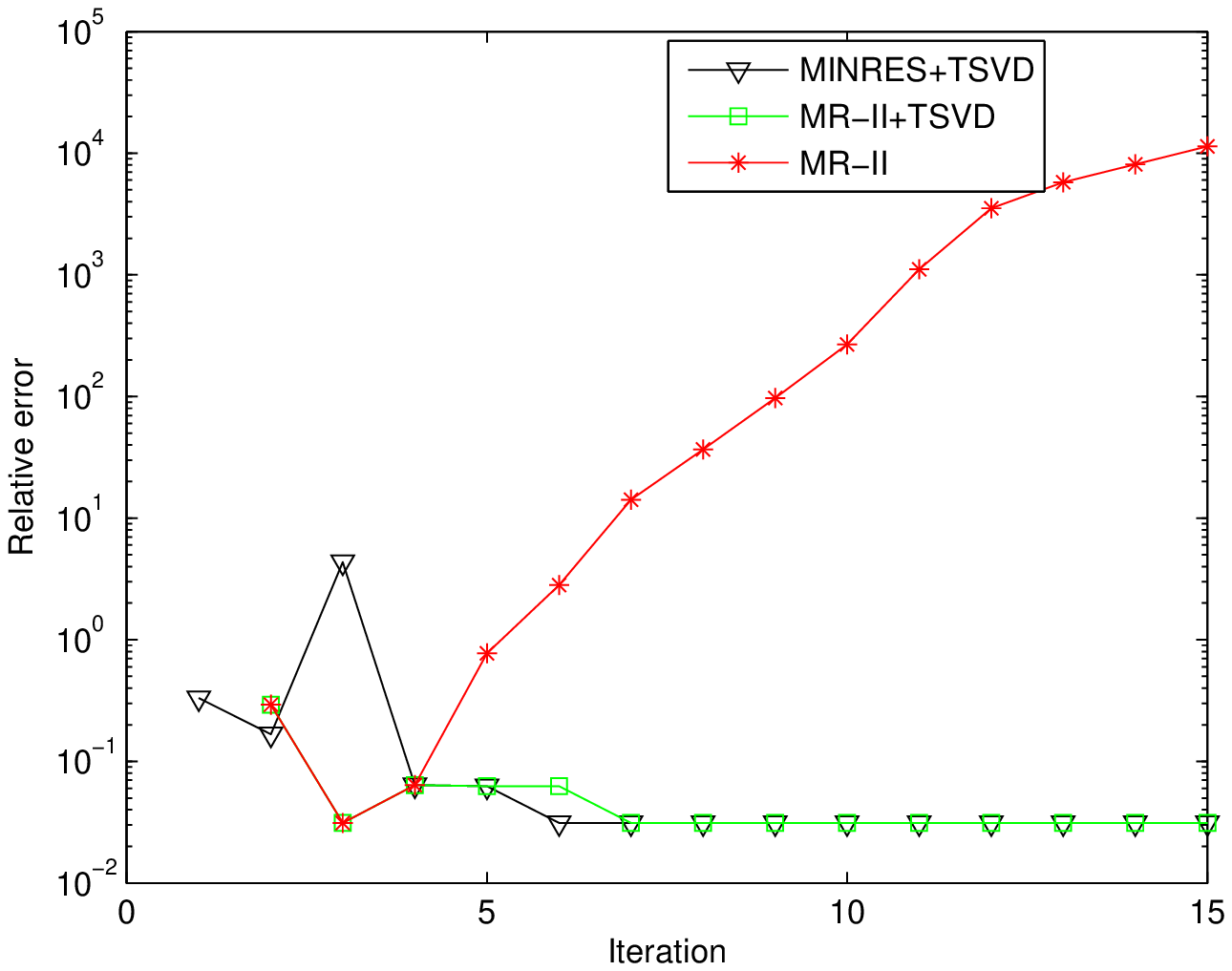}}
  \centerline{(b)}
\end{minipage}
\vfill
\begin{minipage}{0.48\linewidth}
  \centerline{\includegraphics[width=7.0cm,height=5cm]{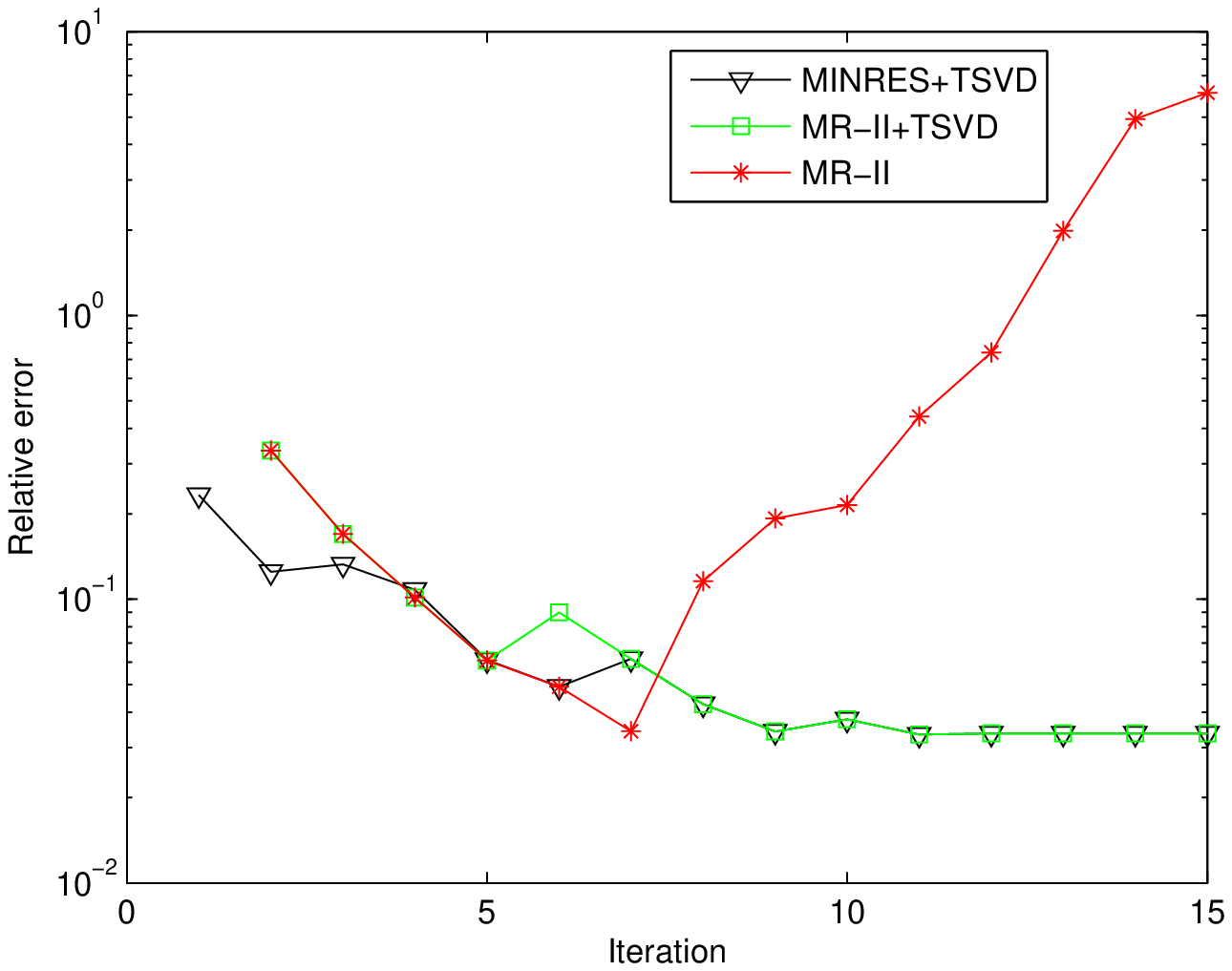}}
  \centerline{(c)}
\end{minipage}
\hfill
\begin{minipage}{0.48\linewidth}
  \centerline{\includegraphics[width=7.0cm,height=5cm]{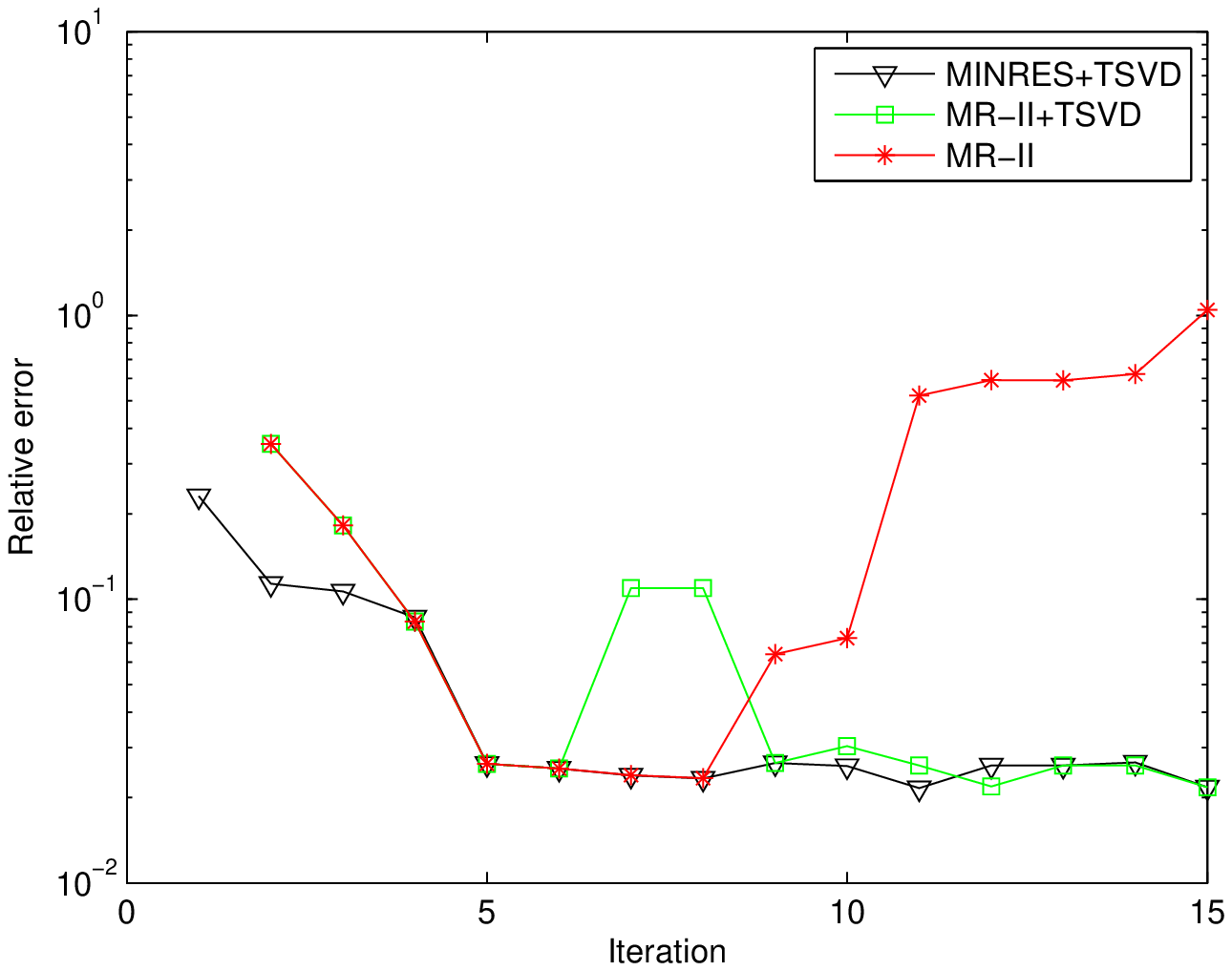}}
  \centerline{(d)}
\end{minipage}
\caption{The relative errors $\|x^{(k)}-x_{true}\|/\|x_{true}\|$
by MR-II, and hybrid MR-II
and MINRES with additional TSVD regularization for
shaw, foxgood, gravity, phillips (from top left to bottom right). }
\label{fig4}
\end{figure}

For MR-II and the hybrid MR-II, we observe from Figure~\ref{fig4}
that MR-II reaches the same error level as the hybrid MR-II,
and the TSVD regularization applied to projected problems simply makes
the regularized solutions with the minimum error almost stabilize
and does not improve the regularized solution by MR-II at semi-convergence.
This justifies the full regularization of MR-II.

Compared with Figures~\ref{fig1}--\ref{fig2},
we find from Figure~\ref{fig4} that the hybrid MINRES improves on MINRES
substantially and the relative errors of iterates
by the hybrid MINRES reach the same minimum level as MR-II and the hybrid MR-II.
These phenomena again justify our
assertion in Section \ref{SectionCom} that the hybrid MINRES is
necessary, independent of the degree of ill-posedness, and
the hybrid MINRES is as effective as the hybrid MR-II.

\begin{figure}

\begin{minipage}{0.48\linewidth}
  \centerline{\includegraphics[width=7.0cm,height=5cm]{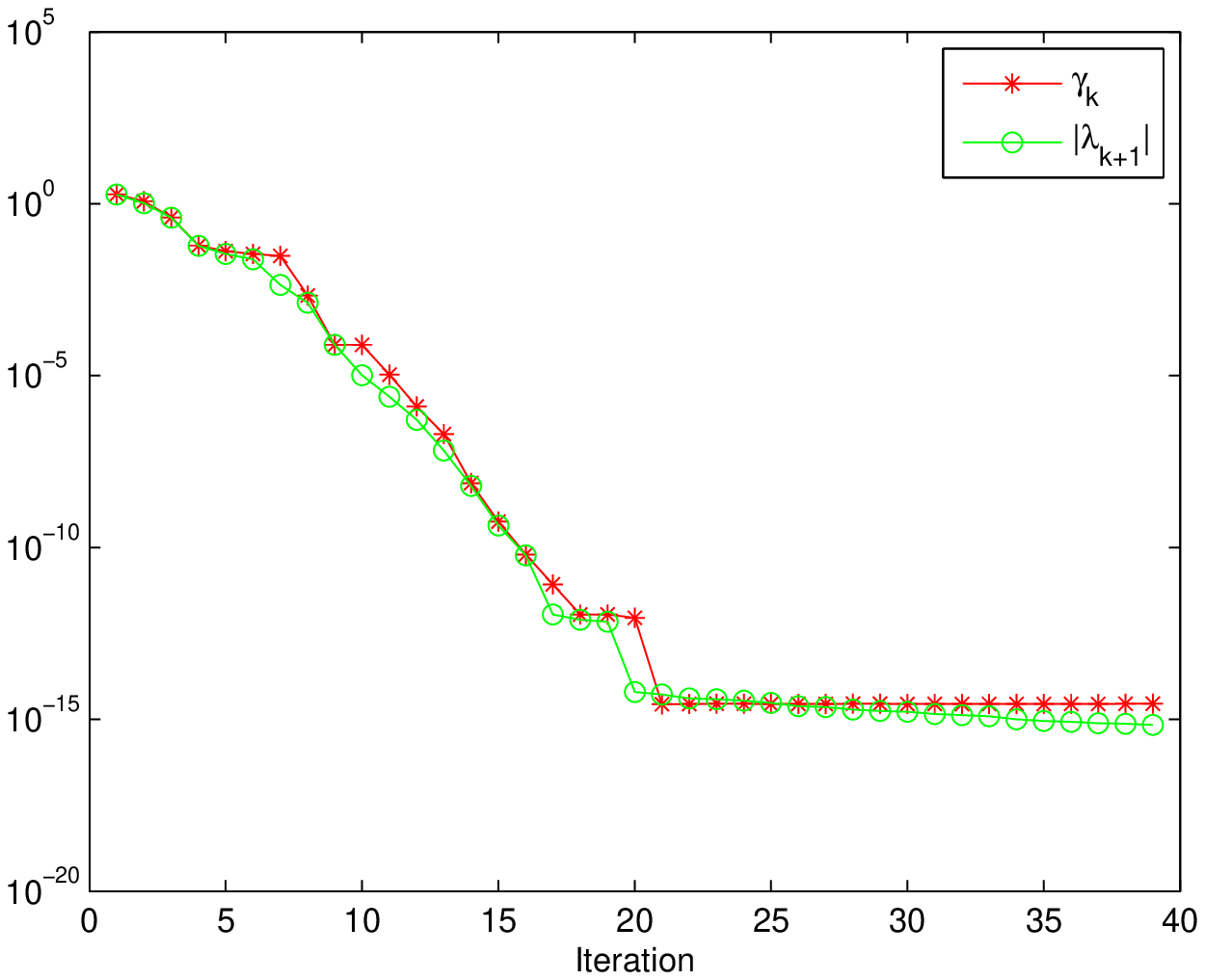}}
  \centerline{(a)}
\end{minipage}
\hfill
\begin{minipage}{0.48\linewidth}
  \centerline{\includegraphics[width=7.0cm,height=5cm]{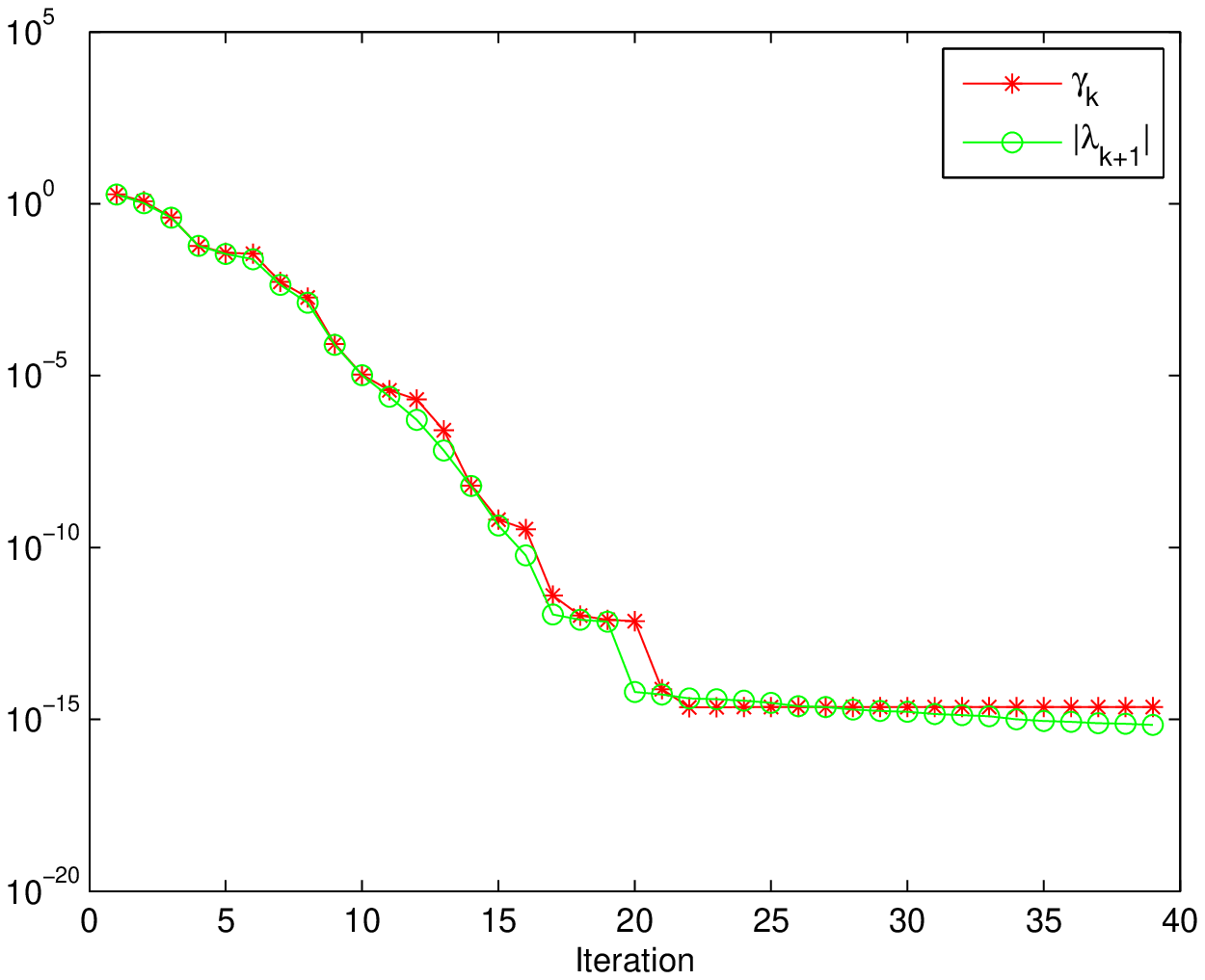}}
  \centerline{(b)}
\end{minipage}
\vfill
\begin{minipage}{0.48\linewidth}
  \centerline{\includegraphics[width=7.0cm,height=5cm]{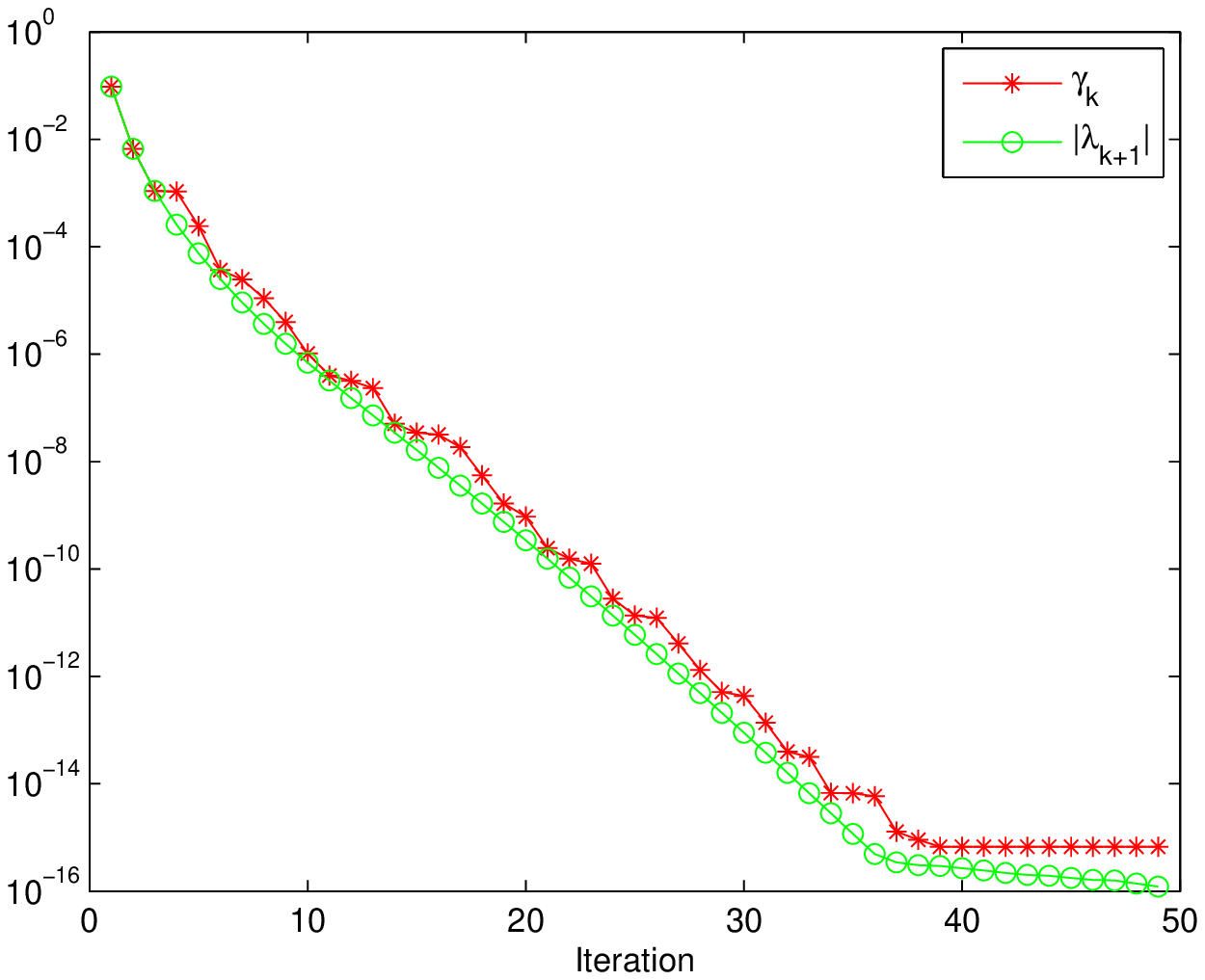}}
  \centerline{(c)}
\end{minipage}
\hfill
\begin{minipage}{0.48\linewidth}
  \centerline{\includegraphics[width=7.0cm,height=5cm]{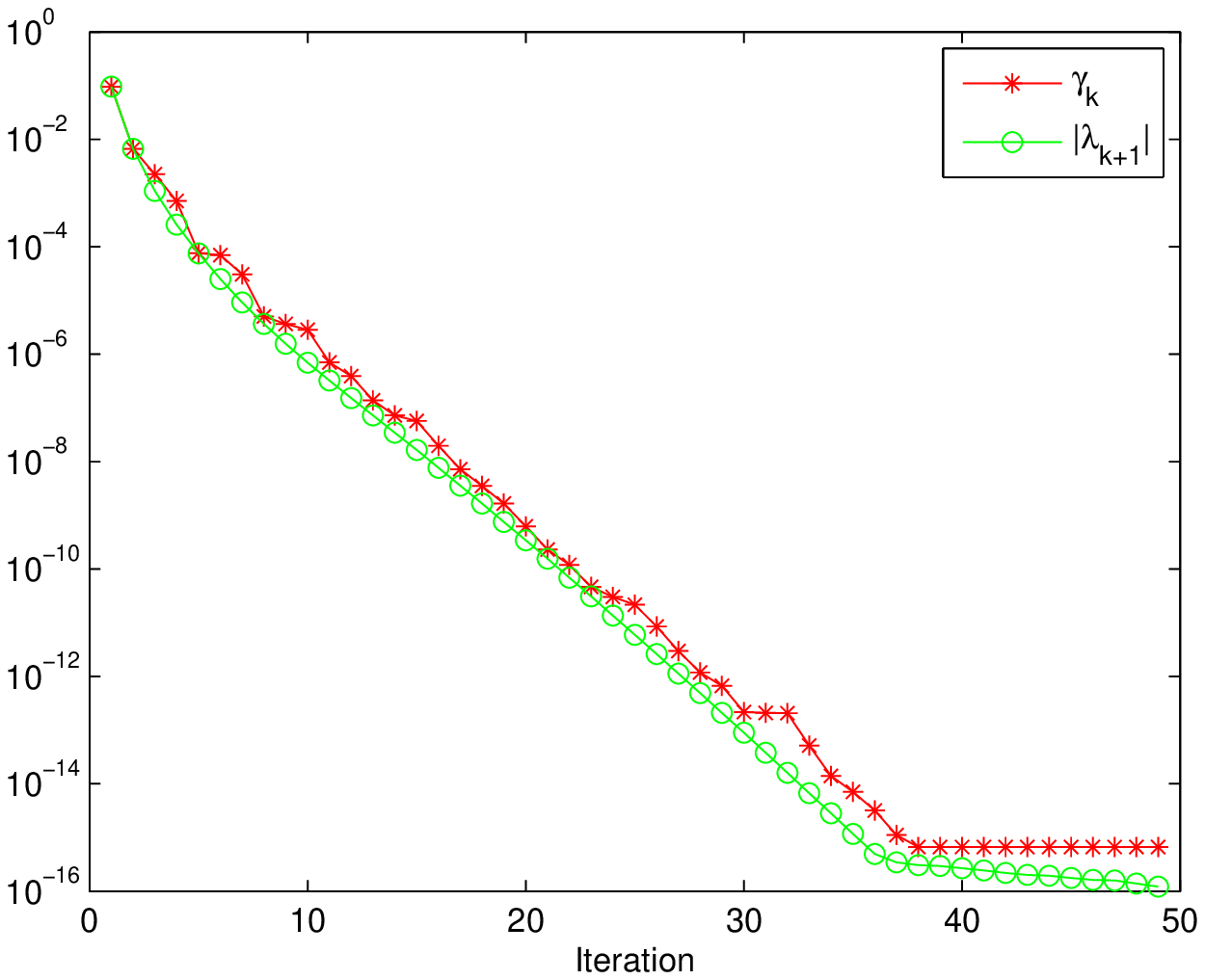}}
  \centerline{(d)}
\end{minipage}
\caption{(a)-(b): Plots of decaying behavior of the sequences $\gamma_k$ and
$|\lambda_{k+1}|$ for shaw with $\varepsilon=10^{-2}$ (left) and
$\varepsilon=10^{-3}$ (right) by MR-II; (c)-(d): Plots of decaying behavior of the
sequences $\gamma_k$ and $|\lambda_{k+1}|$ for foxgood with
$\varepsilon=10^{-3}$ (left) and $\varepsilon=10^{-4}$ (right) by MR-II.}
\label{fig5}
\end{figure}

Figure~\ref{fig5} and Figure~\ref{fig7} display the
curves of sequences
$\gamma_k$ with the noise levels $\varepsilon=10^{-2}, 10^{-3}, 10^{-4}$,
respectively, for the four severely and moderately problems.
We see that $\gamma_k\approx\sigma_{k+1}=|\lambda_{k+1}|$,
almost independent of noise level $\varepsilon$. We point out that,
due to the round-offs in finite precision arithmetic, they level off at
the level of $\epsilon_{\rm mach}$ when $k=20$ for shaw, $k=37$ for foxgood and
$k=50$ for gravity. The results indicate that the $Q_{k+1}T_kQ_k^T$
are near best rank $k$ approximations to $A$ with the approximate
accuracy $\sigma_{k+1}$ so that $T_k$ does not become ill-conditioned
before $k\leq k_0$. As a result, the regularized solutions $x^{(k)}$
become increasingly better approximations to $x_{true}$ until iteration $k_0$,
and they are deteriorated after that iteration.
At iteration $k_0$, $x^{(k_0)}$ captures the $k_0$ dominant spectral or
equivalent SVD components of $A$ and is a best possible regularized solution,
i.e., MR-II has the full regularization for the severely
ill-posed problems tested.

\begin{figure}
\begin{minipage}{0.48\linewidth}
  \centerline{\includegraphics[width=7.0cm,height=5cm]{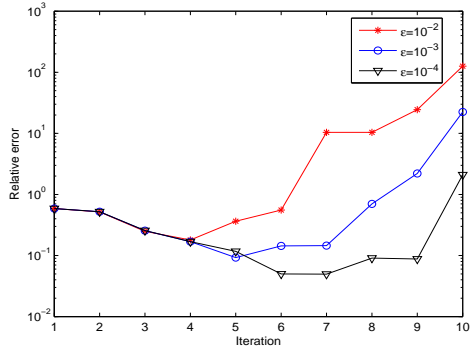}}
  \centerline{(a)}
\end{minipage}
\hfill
\begin{minipage}{0.48\linewidth}
  \centerline{\includegraphics[width=7.0cm,height=5cm]{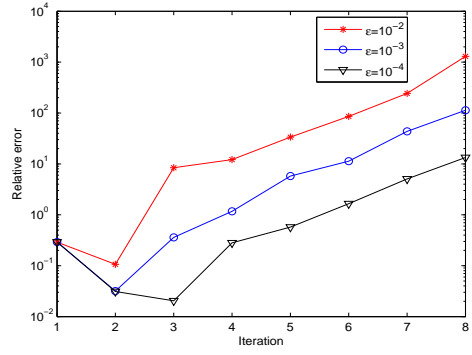}}
  \centerline{(b)}
\end{minipage}
\caption{ The relative error $\left\|x^{(k)}-x_{true}\right\|/
\|x_{true}\|$ with respect to $\varepsilon=10^{-2}, 10^{-3}, 10^{-4}$
for shaw (left) and foxgood (right) by MR-II.}
\label{fig6}
\end{figure}

\begin{figure}
\begin{minipage}{0.48\linewidth}
  \centerline{\includegraphics[width=7.0cm,height=5cm]{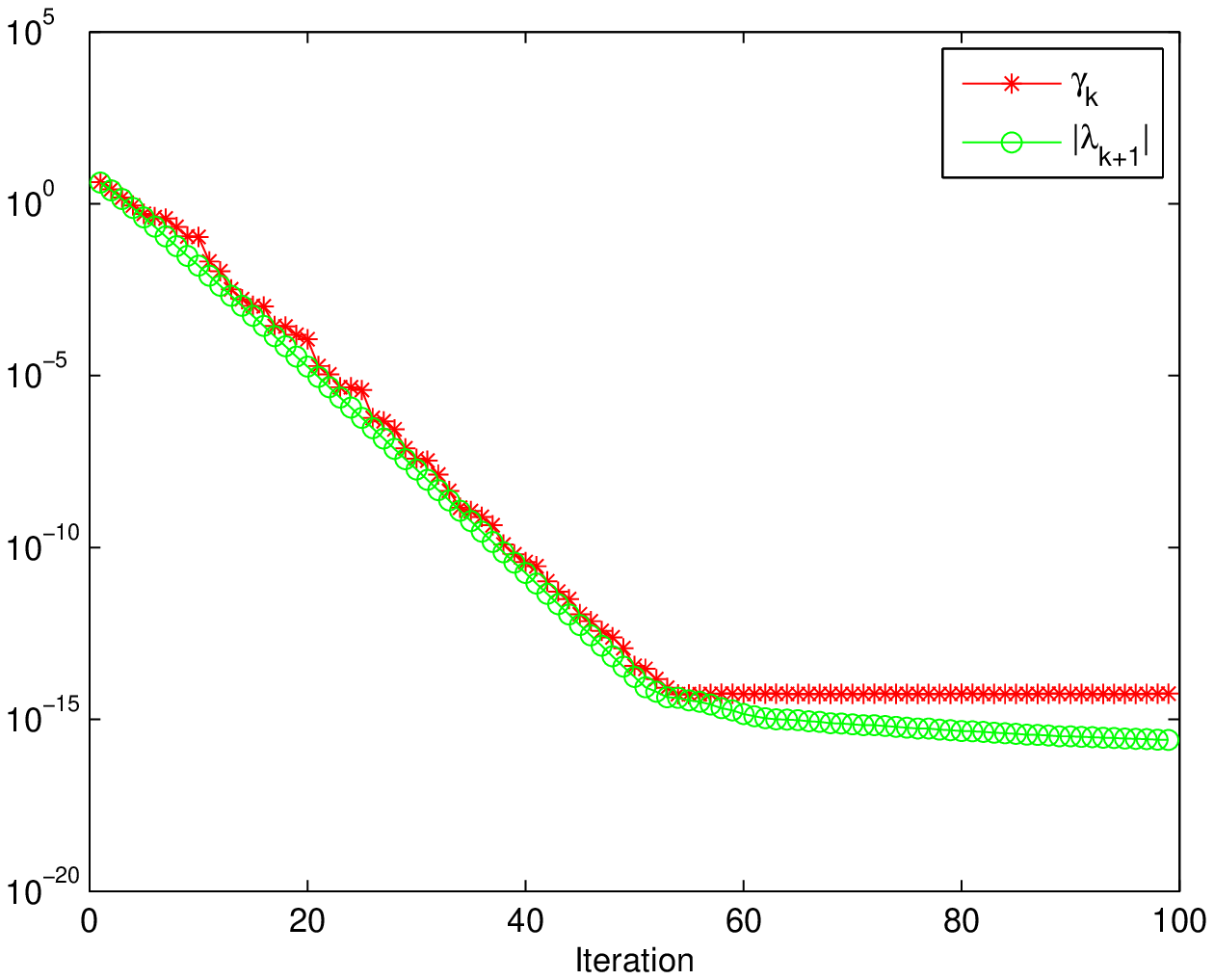}}
  \centerline{(a)}
\end{minipage}
\hfill
\begin{minipage}{0.48\linewidth}
  \centerline{\includegraphics[width=7.0cm,height=5cm]{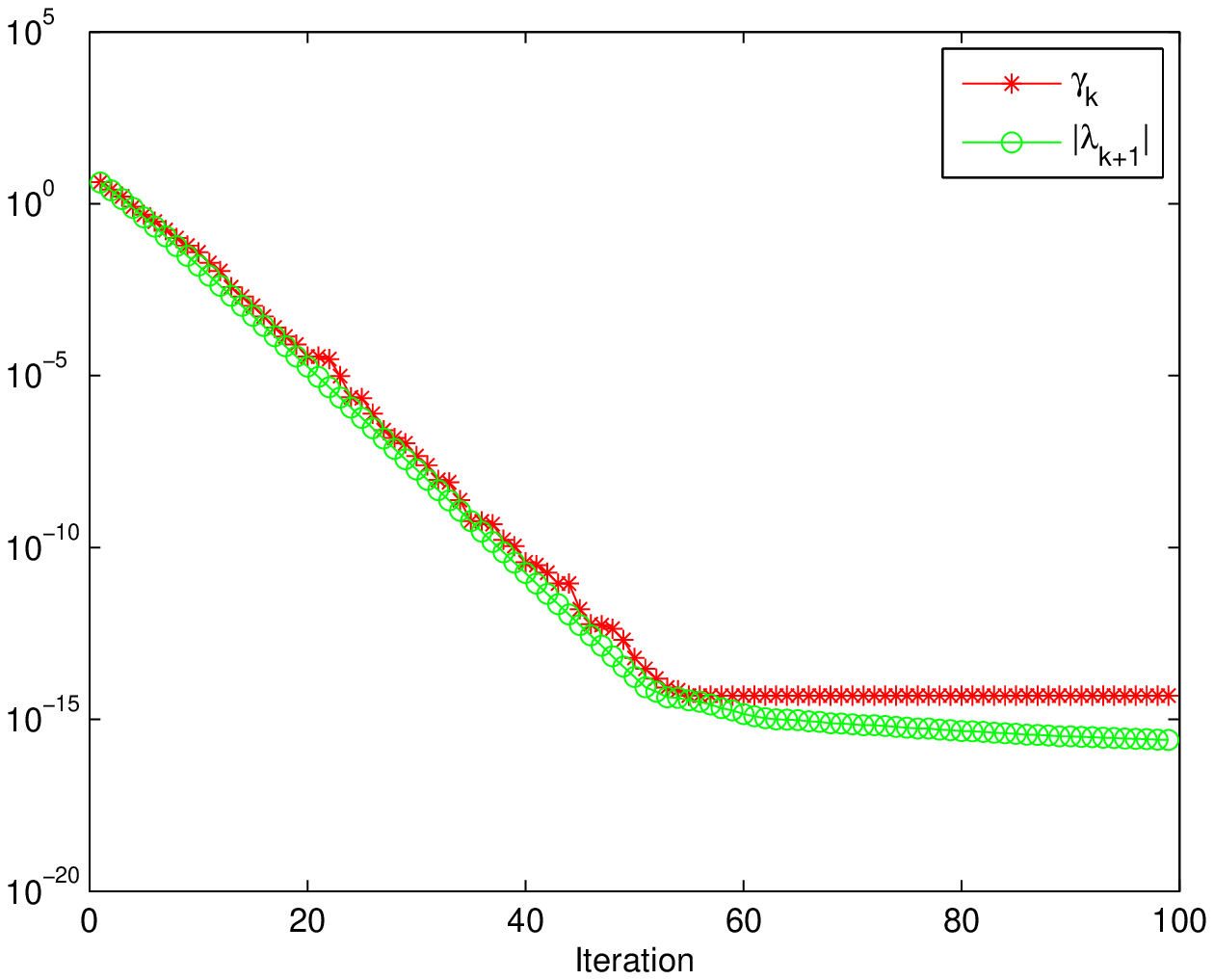}}
  \centerline{(b)}
\end{minipage}
\vfill
\begin{minipage}{0.48\linewidth}
  \centerline{\includegraphics[width=7.0cm,height=5cm]{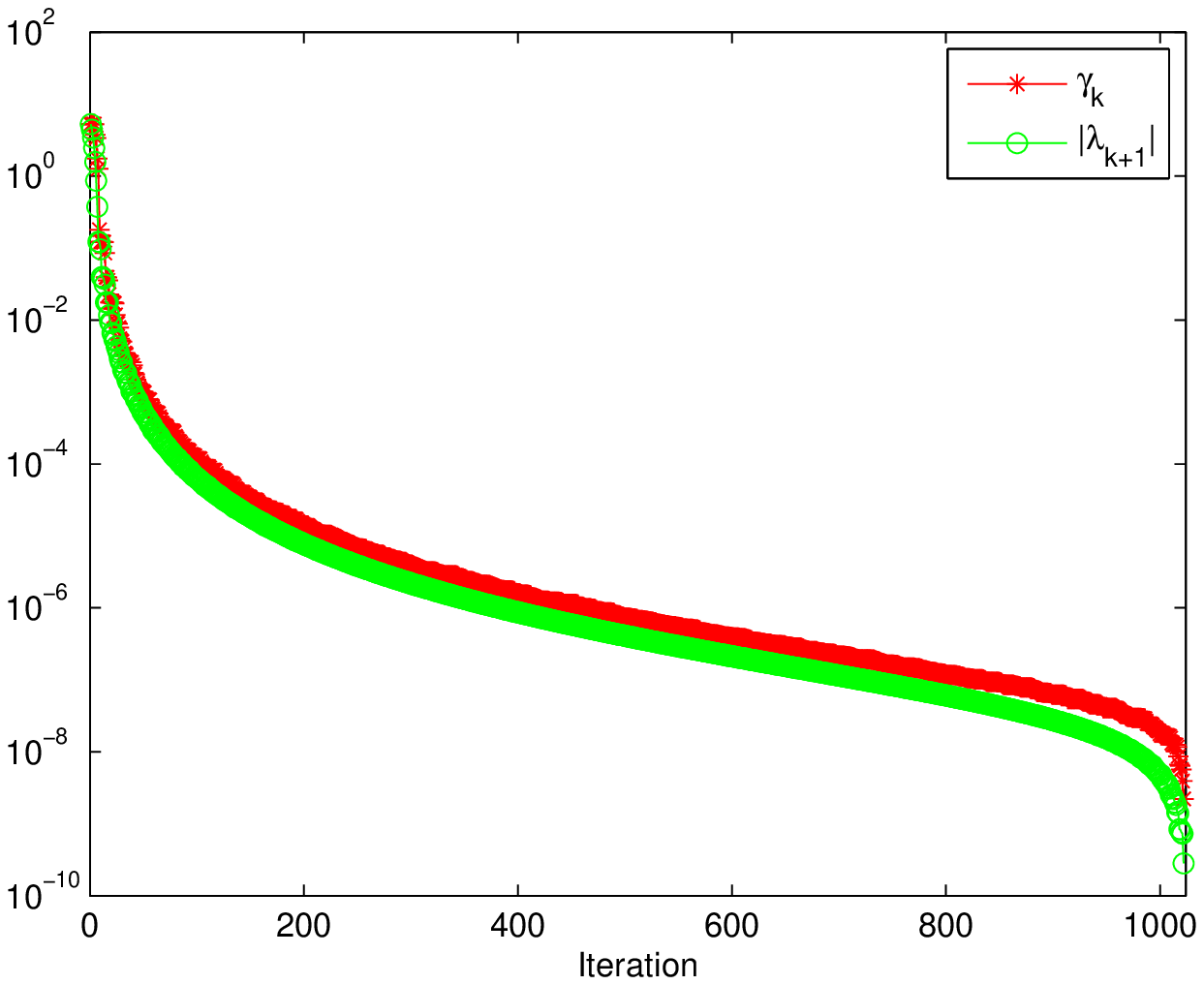}}
  \centerline{(c)}
\end{minipage}
\hfill
\begin{minipage}{0.48\linewidth}
  \centerline{\includegraphics[width=7.0cm,height=5cm]{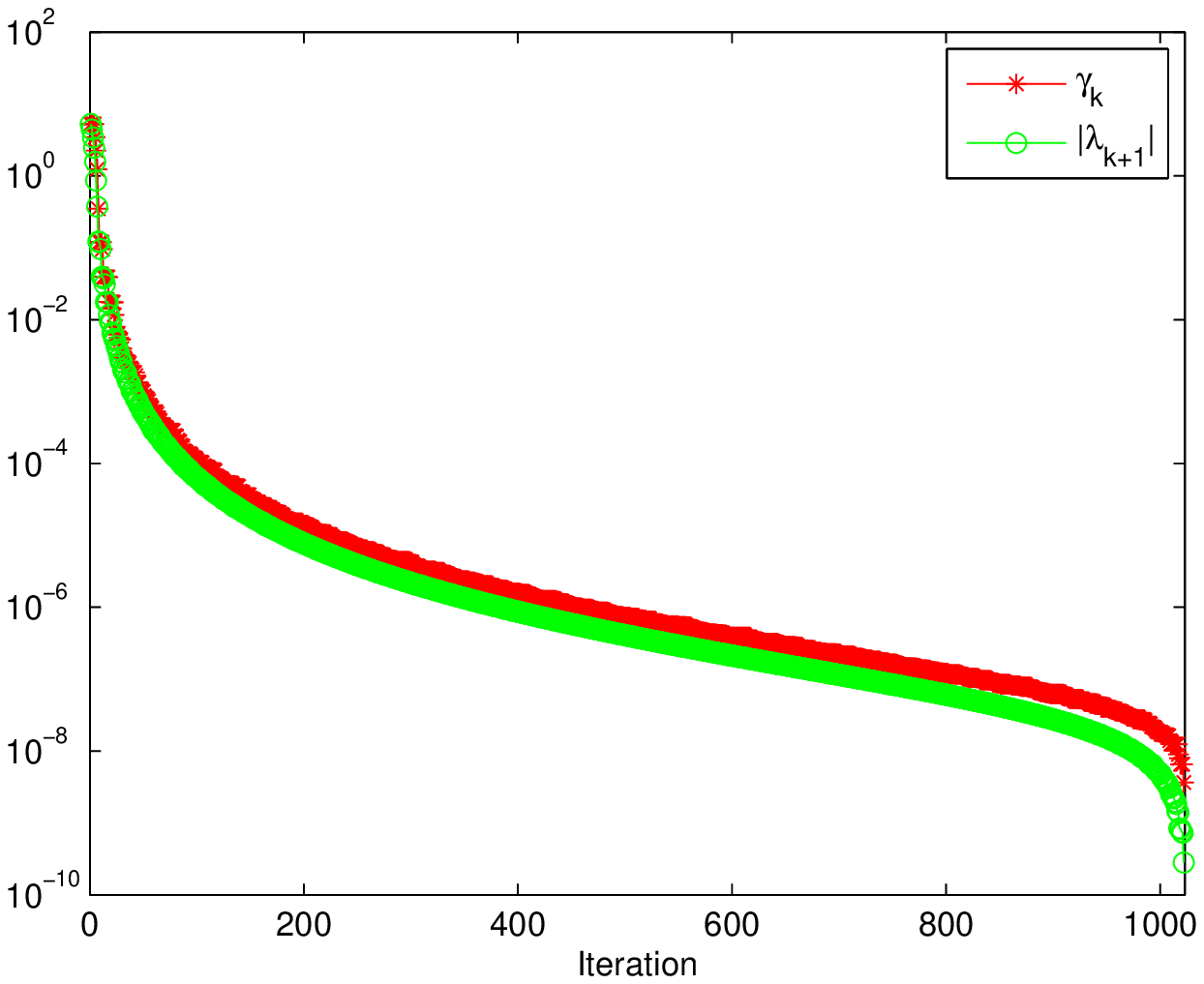}}
  \centerline{(d)}
\end{minipage}
\caption{(a)-(b): Plots of decaying behavior of the sequences $\gamma_k$
and $|\lambda_{k+1}|$ for gravity with $\varepsilon=10^{-2}$ (left)
and $\varepsilon=10^{-3}$ (right) by MR-II; (c)-(d): Plots of decaying behavior of
the sequences $\gamma_k$ and $|\lambda_{k+1}|$ for phillips
with $\varepsilon=10^{-3}$ (left) and $\varepsilon=10^{-4}$ (right) by MR-II.}
\label{fig7}
\end{figure}

\begin{figure}
\begin{minipage}{0.48\linewidth}
  \centerline{\includegraphics[width=7.0cm,height=5cm]{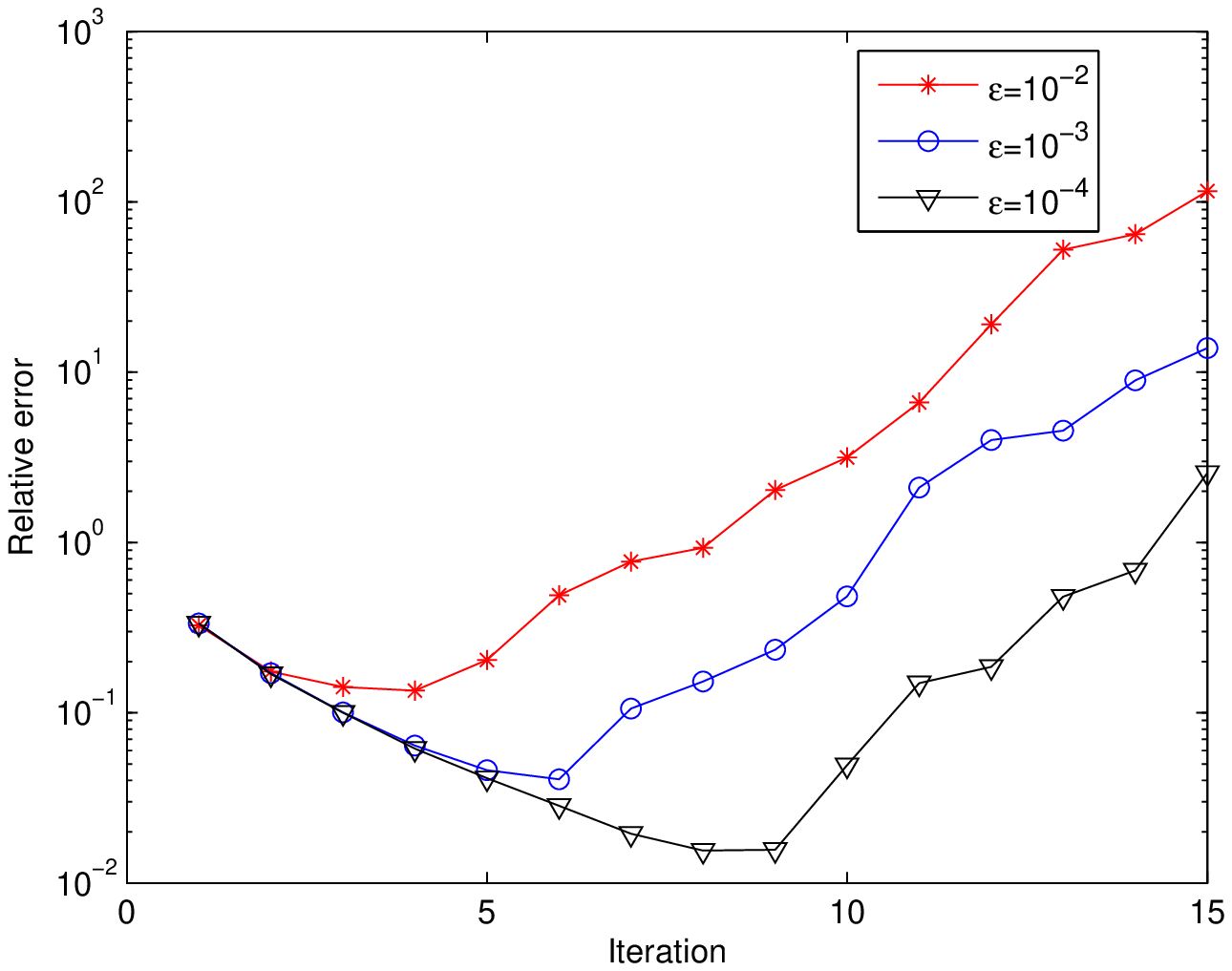}}
  \centerline{(a)}
\end{minipage}
\hfill
\begin{minipage}{0.48\linewidth}
  \centerline{\includegraphics[width=7.0cm,height=5cm]{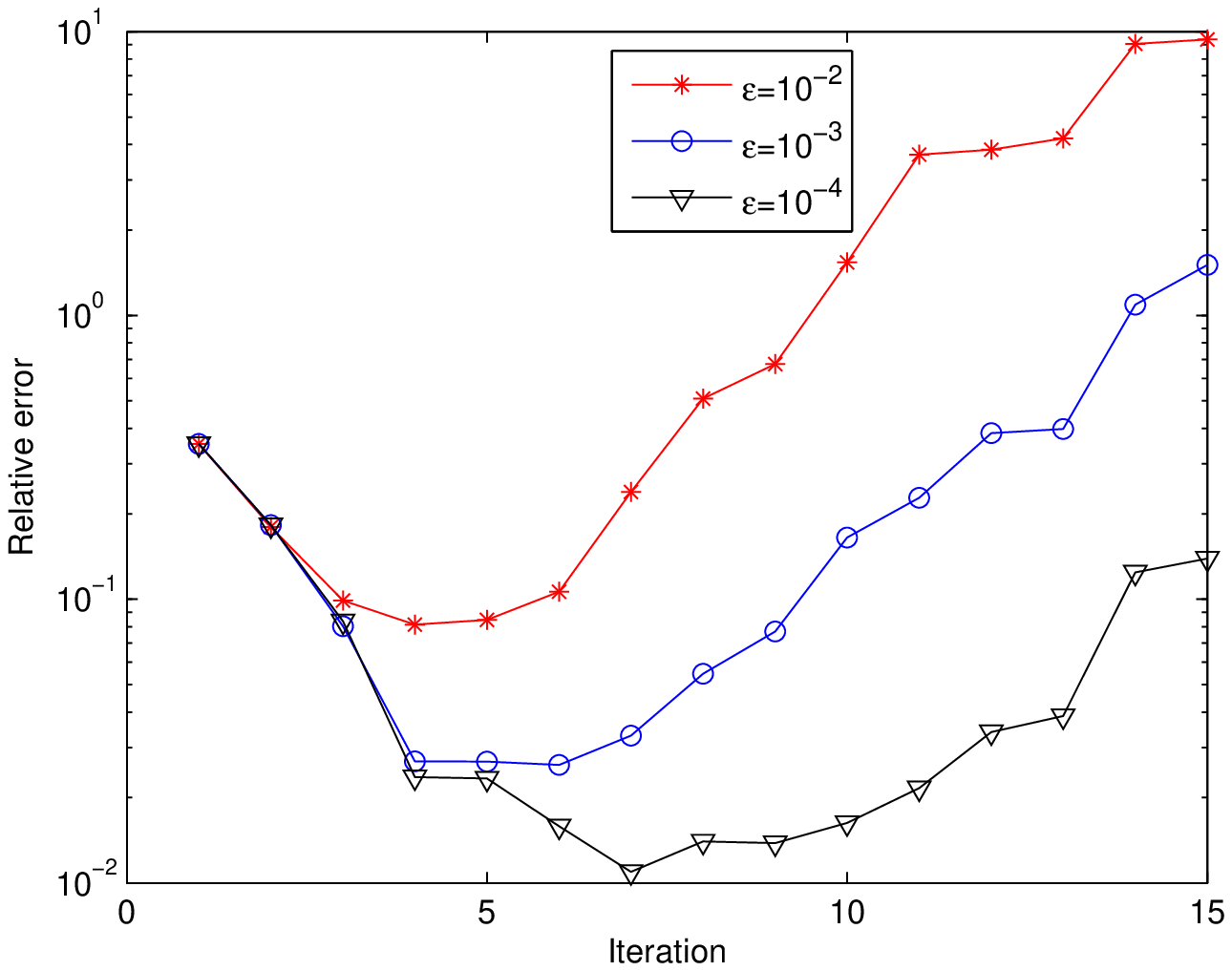}}
  \centerline{(b)}
\end{minipage}
\caption{ The relative errors $\left\|x^{(k)}-x_{true}\right\|/
\|x_{true}\|$ with respect to $\varepsilon=10^{-2}, 10^{-3}, 10^{-4}$
for gravity (left) and phillips (right) by MR-II.}
\label{fig8}
\end{figure}

Figure~\ref{fig6} and Figure~\ref{fig8} plot the relative errors
$\left\|x^{(k)}-x_{true}\right\|/\|x_{true}\|$ with different noise levels
for these four severely and moderately ill-posed problems. For smaller noise
levels, MR-II gets more accurate best
regularized solutions at cost of more iterations. This is expected since, from
\eqref{picard} and $|\lambda_{k_0+1}^{1+\beta}|=\mid v_{k_0+1}^T \hat b\mid
\approx \mid v_{k_0+1}^T \hat b\mid$,
a bigger $k_0$ is needed for a smaller noise level $\|e\|$. Moreover, MR-II
needs more iterations to achieve semi-convergence for moderately ill-posed
problems with the same noise level, since $\sigma_j$ does not decay as fast
as that for a severely ill-posed problem.

\begin{figure}[t]
\begin{minipage}{0.48\linewidth}
  \centerline{\includegraphics[width=7.0cm,height=5cm]{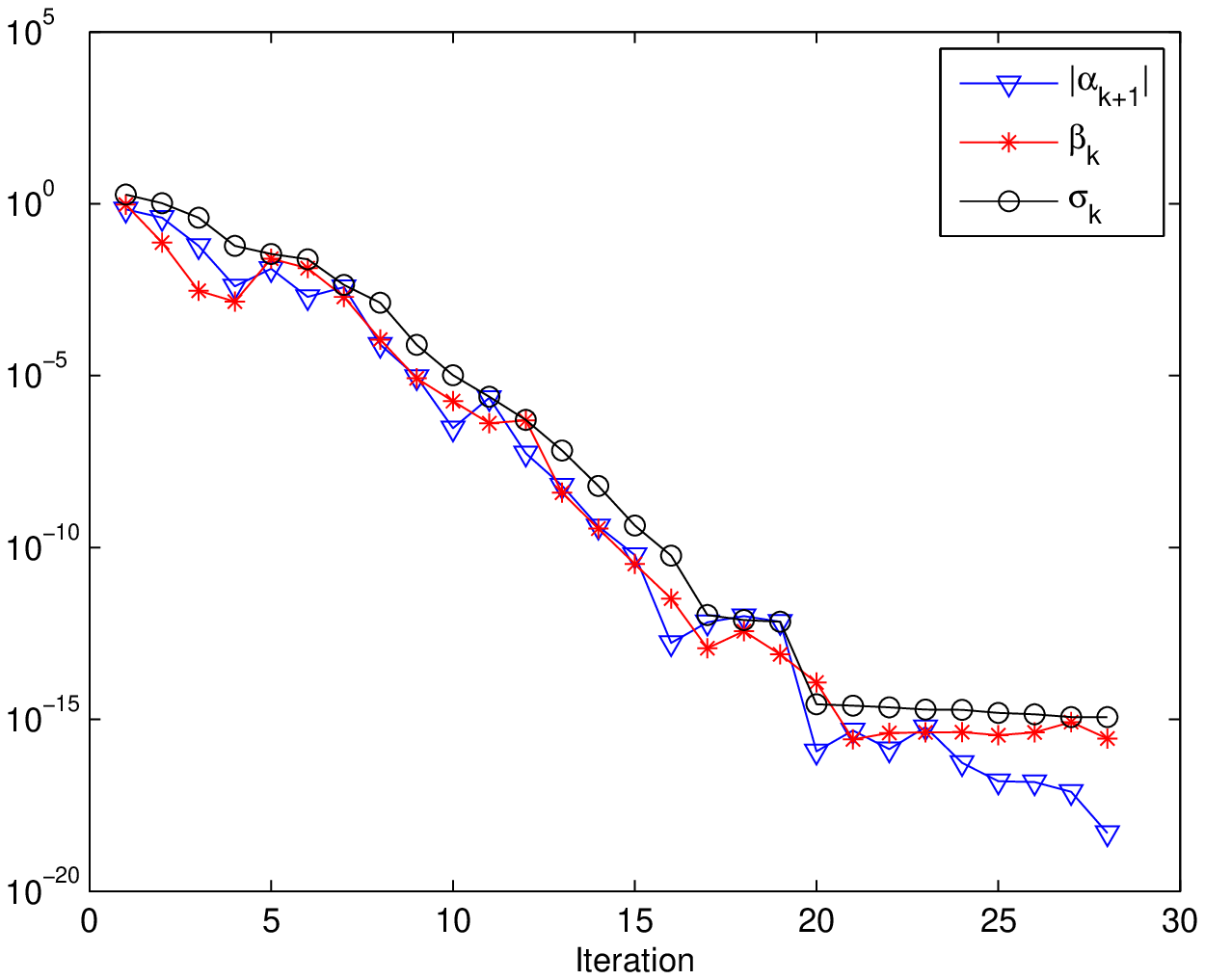}}
  \centerline{(a)}
\end{minipage}
\hfill
\begin{minipage}{0.48\linewidth}
  \centerline{\includegraphics[width=7.0cm,height=5cm]{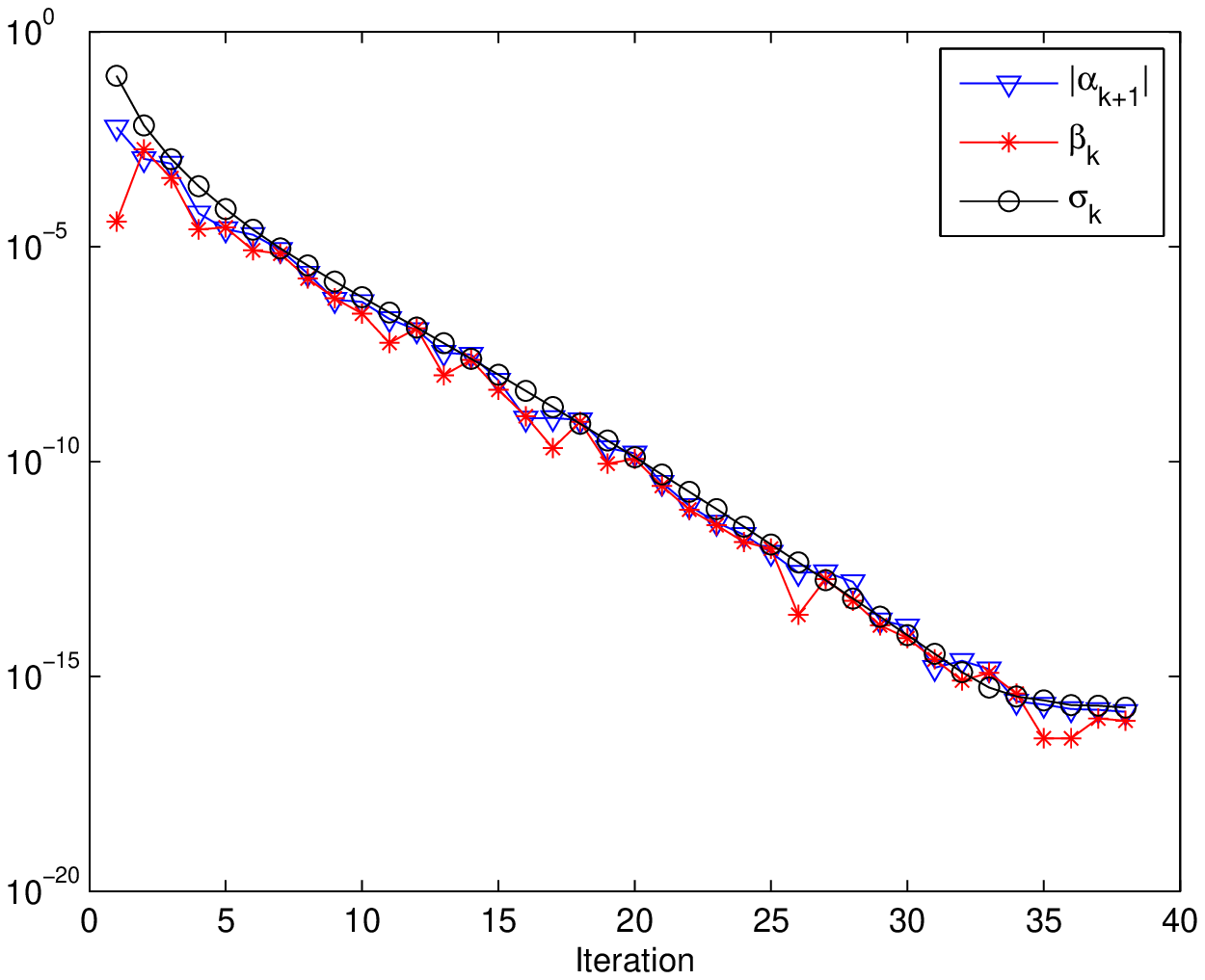}}
  \centerline{(b)}
\end{minipage}
\vfill
\begin{minipage}{0.48\linewidth}
  \centerline{\includegraphics[width=7.0cm,height=5cm]{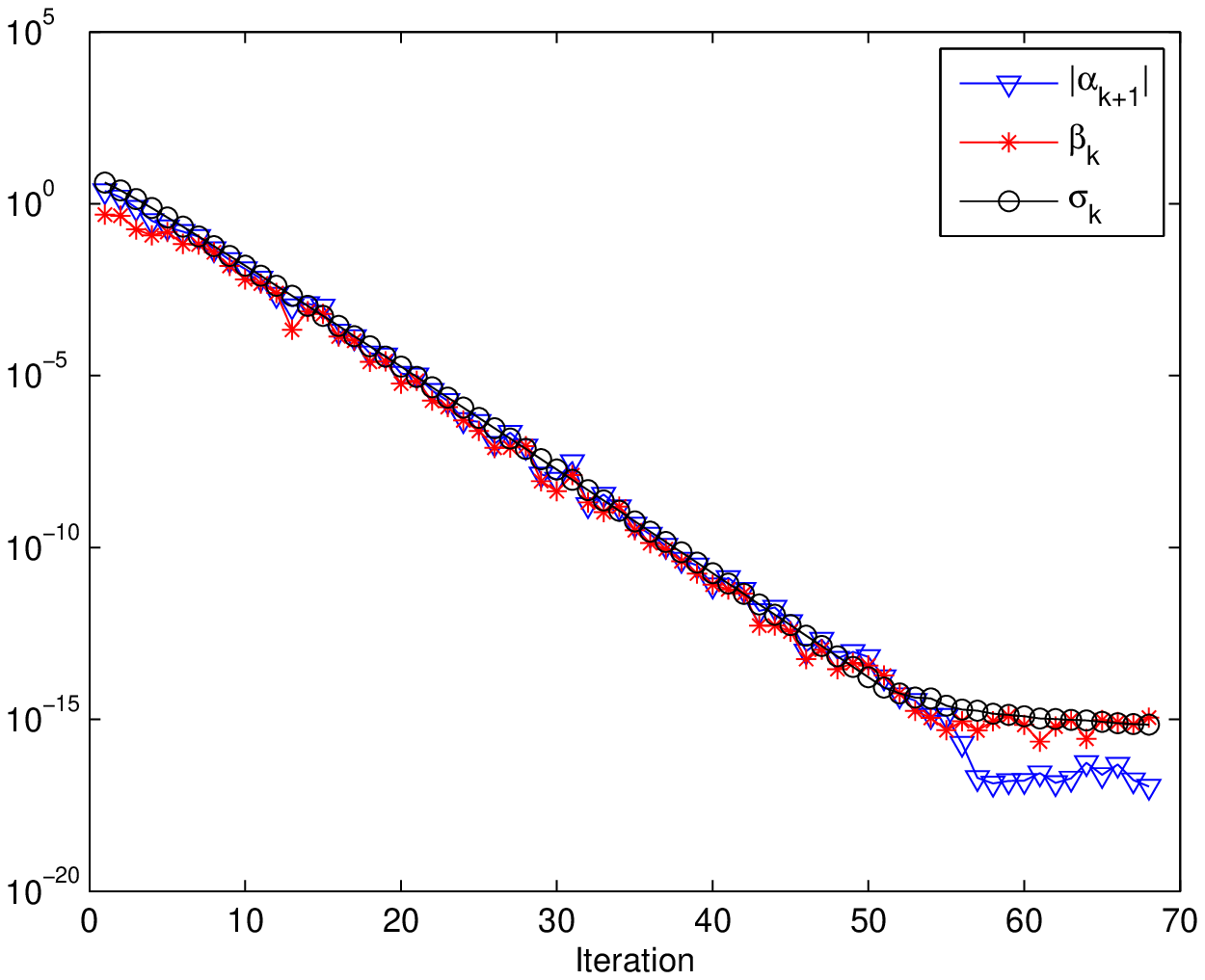}}
  \centerline{(c)}
\end{minipage}
\hfill
\begin{minipage}{0.48\linewidth}
  \centerline{\includegraphics[width=7.0cm,height=5cm]{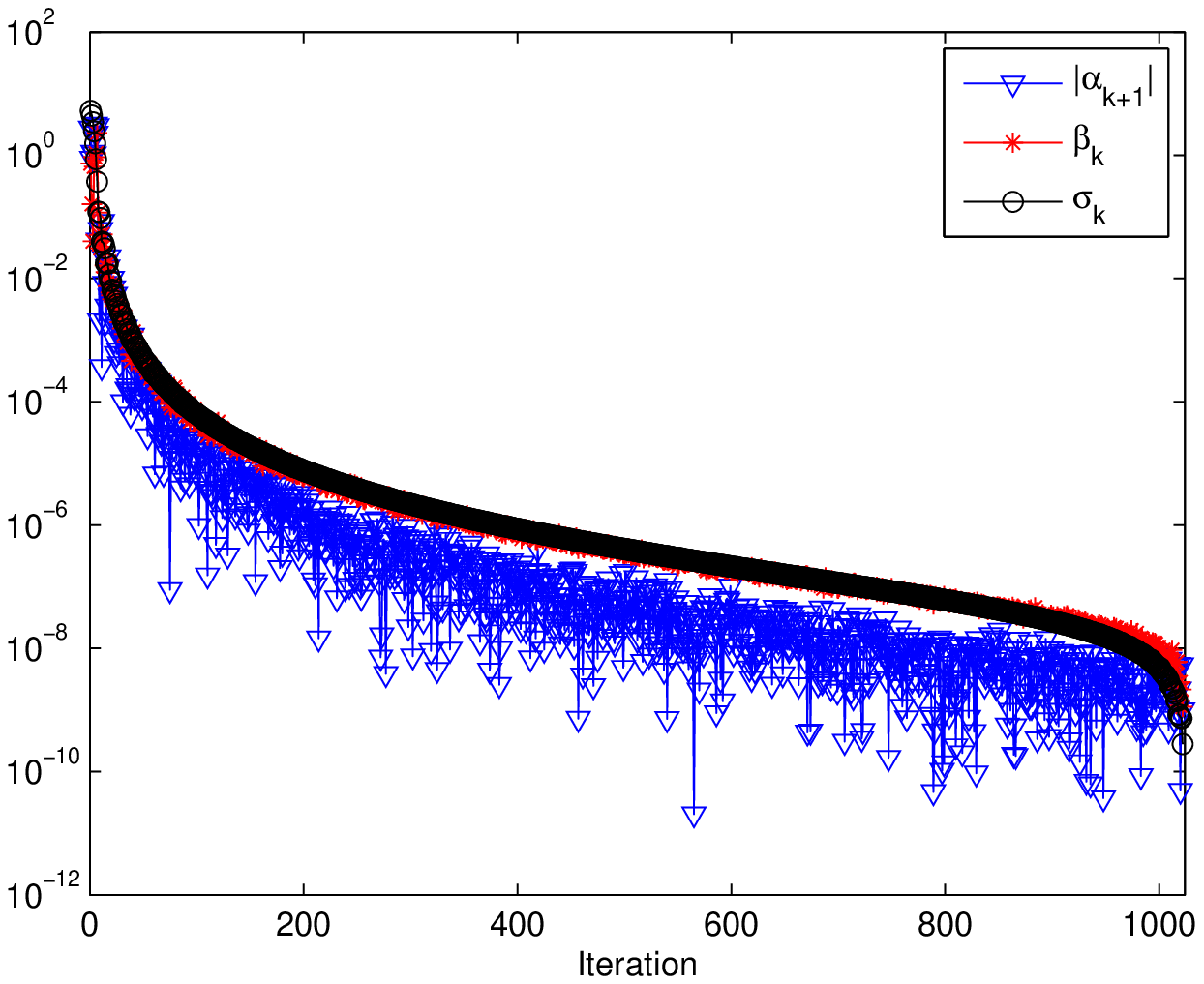}}
  \centerline{(d)}
\end{minipage}
\caption{(a)-(d): Plots of decaying behavior of the sequences $|\alpha_{k+1}|$, $\beta_k$
and $\sigma_{k}$ for shaw, foxgood, gravity, phillips (from top left to bottom right)
with $\varepsilon=10^{-3}$ by MR-II.}
\label{figpl}
\end{figure}

Figures~\ref{figpl} display the decreasing curves of quantities $|\alpha_{k+1}|$,
$\beta_k$ and $\sigma_k$, $k=2,\ldots,n-1$. From Figure~\ref{figpl} (a),
we see that, for the severely ill-posed problem shaw,
$|\alpha_{k+1}|$ and $\beta_k$ decrease as fast as $\sigma_{k}$
and the three quantities level off at the level of $\epsilon_{\rm mach}$ for $k$ no
more than 20, and after that these quantities are purely round-offs and not
reliable any more. Similar phenomena are also observed for the other two
severely ill-posed problems foxgood and gravity, as indicated
by Figure~\ref{figpl} (b) and (c). Figure~\ref{figpl} (d) illustrates
that $\beta_k$ decreases
as fast as $\sigma_{k}$ but $|\alpha_{k+1}|$ decays as fast as $\sigma_k$
in the first iterations and then considerably faster than $\sigma_k$ as
$k$ increases in the later stage for moderately ill-posed problem phillips.

Finally, we report some comparison results on LSQR, the hybrid LSQR and
MR-II, the hybrid MR-II.  As already proved in \cite{huang14},
a hybrid LSQR should be used to compute best possible regularized
solutions for mildly ill-posed problems. It has also been experimentally justified
in \cite{huang14} that LSQR has the full regularization for severely and moderately
ill-posed problems. We have tested LSQR and the hybrid LSQR, and compared their
effectiveness and efficiency with MR-II and the hybrid MR-II.
We have found that, for each of the severely and moderately ill-posed problems
in Table~\ref{table} and with the same noise level,
both the pure MR-II and LSQR obtain the best regularized solutions with almost
the same accuracy using almost the same iterations.
For deriv2, the hybrid MR-II and LSQR compute the best possible
regularized solutions using almost the same iterations. These results
tell us two things: (i) As an iterative regularization method,
MR-II is as effective as LSQR for an ill-posed
problem; (ii) MR-II is twice as efficient as LSQR.

\subsection{A 2D image restoration problem}

The problem blur is a 2D image deblurring problem and more complex than
the other five 1D problems in Table~\ref{table}. It arises in connection with
the degradation of digital images by atmospheric turbulence blur.
We use the code $\mathsf{blur(n,band,sigma)}$
in \cite{hansen07} to generate an $n^2\times n^2$ matrix $A$,
the true solution $x_{true}$ and noise-free right-hand $\hat{b}$.
The vector $x_{true}$ is a columnwise stacked
version of a simple test image, while $\hat{b}=Ax_{true}$
holds for a columnwise stacked version of the blurred image.
The blurring matrix $A$ is block Toeplitz with Toeplitz blocks, which
has two parameters $\mathsf{band}$ and $\mathsf{sigma}$;
the former specifies the half-bandwidth of the Toeplitz blocks,
and the latter controls the shape of the Gaussian
point spread function and thus the amount of smoothing.
We generate a blurred and noisy image $b=\hat{b}+e$ by
adding a while noise vector $e$. The goal is to restore
the true image $x_{true}$ from $b$.

We take $n=256$ and the relative noise level $\varepsilon=5\times10^{-3}$,
giving rise to $A$ with order $n^2=65,536$.
It is known that the larger the $\mathsf{sigma}$,
the less ill-posed the problem. Purely for an experimental
purpose, we computed all the singular values of a few $A$
with $n^2\leq 10,000$ using the matlab function {\sf svd}.
Since the degree of ill-posedness is the same for different large
$n^2$, we have deduced from the computed singular
values for these matrices $A$ that
$\mathsf{band=3}, \mathsf{sigma=0.7}$ (the default setting)
generates mildly ill-posed problems, while
$\mathsf{band=7}, \mathsf{sigma=2}$ gives rise to moderately ill-posed problems.
We next test MINRES, MR-II and their hybrid variants for these two problems,
and verify the regularizing effects similar to the previous
mildly and moderately ill-posed problems: deriv2 and phillips.

 \begin{figure}[t]
\begin{minipage}{0.48\linewidth}
  \centerline{\includegraphics[width=6.0cm,height=4.5cm]{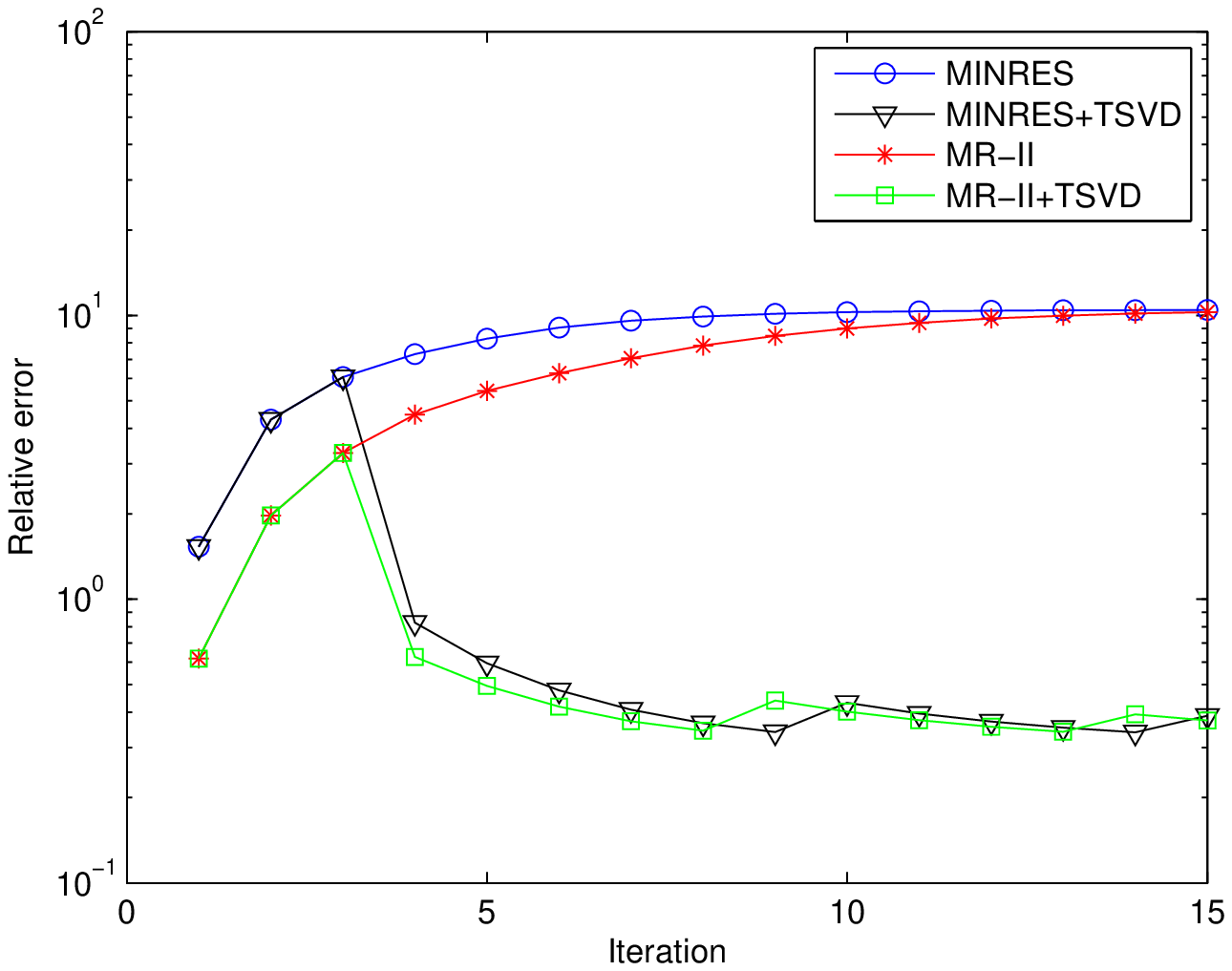}}
  \centerline{(a)}
\end{minipage}
\hfill
\begin{minipage}{0.48\linewidth}
  \centerline{\includegraphics[width=7.0cm,height=5cm]{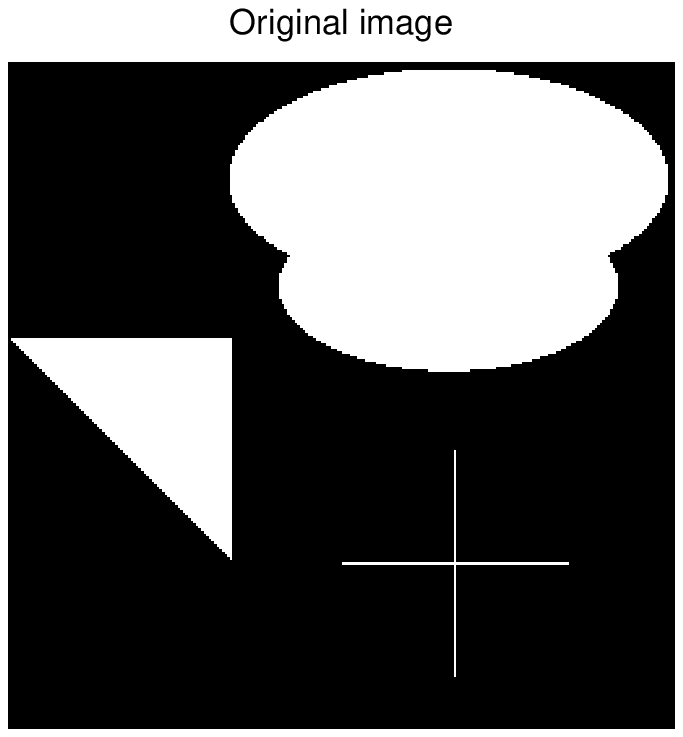}}
  \centerline{(b)}
\end{minipage}
\vfill
\begin{minipage}{0.48\linewidth}
  \centerline{\includegraphics[width=7.0cm,height=5cm]{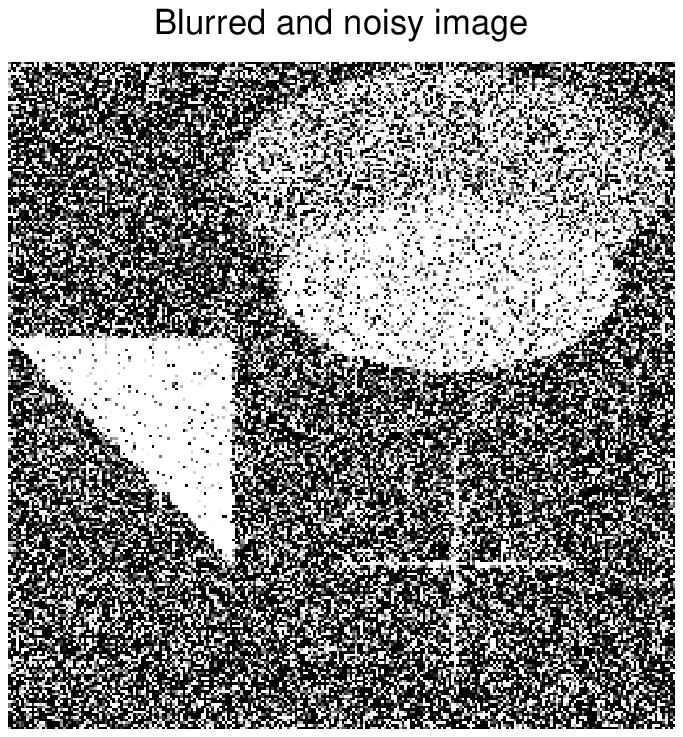}}
  \centerline{(c)}
\end{minipage}
\hfill
\begin{minipage}{0.48\linewidth}
  \centerline{\includegraphics[width=7.0cm,height=5cm]{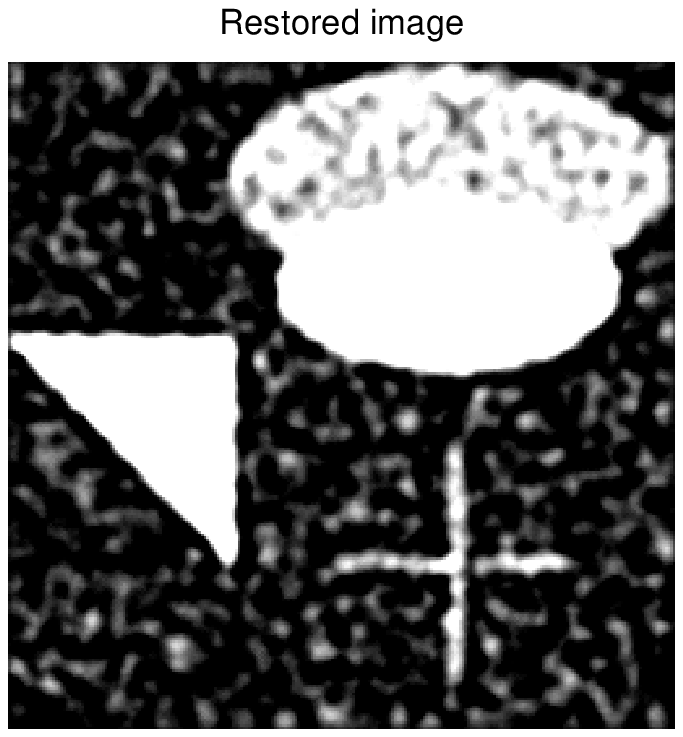}}
  \centerline{(d)}
\end{minipage}
\caption{(a): The relative errors $\|x^{(k)}-x_{true}\|/\|x_{true}\|$
by MINRES, hybrid MINRES, MR-II,
and hybrid MR-II; (b): The original image; (c): The blurred and noisy image;
(d): The restored image with $\mathsf{band=3}, \mathsf{sigma=0.7}$.}
\label{fig11}
\end{figure}

\begin{figure}[t]
\begin{minipage}{0.48\linewidth}
  \centerline{\includegraphics[width=6.0cm,height=4.5cm]{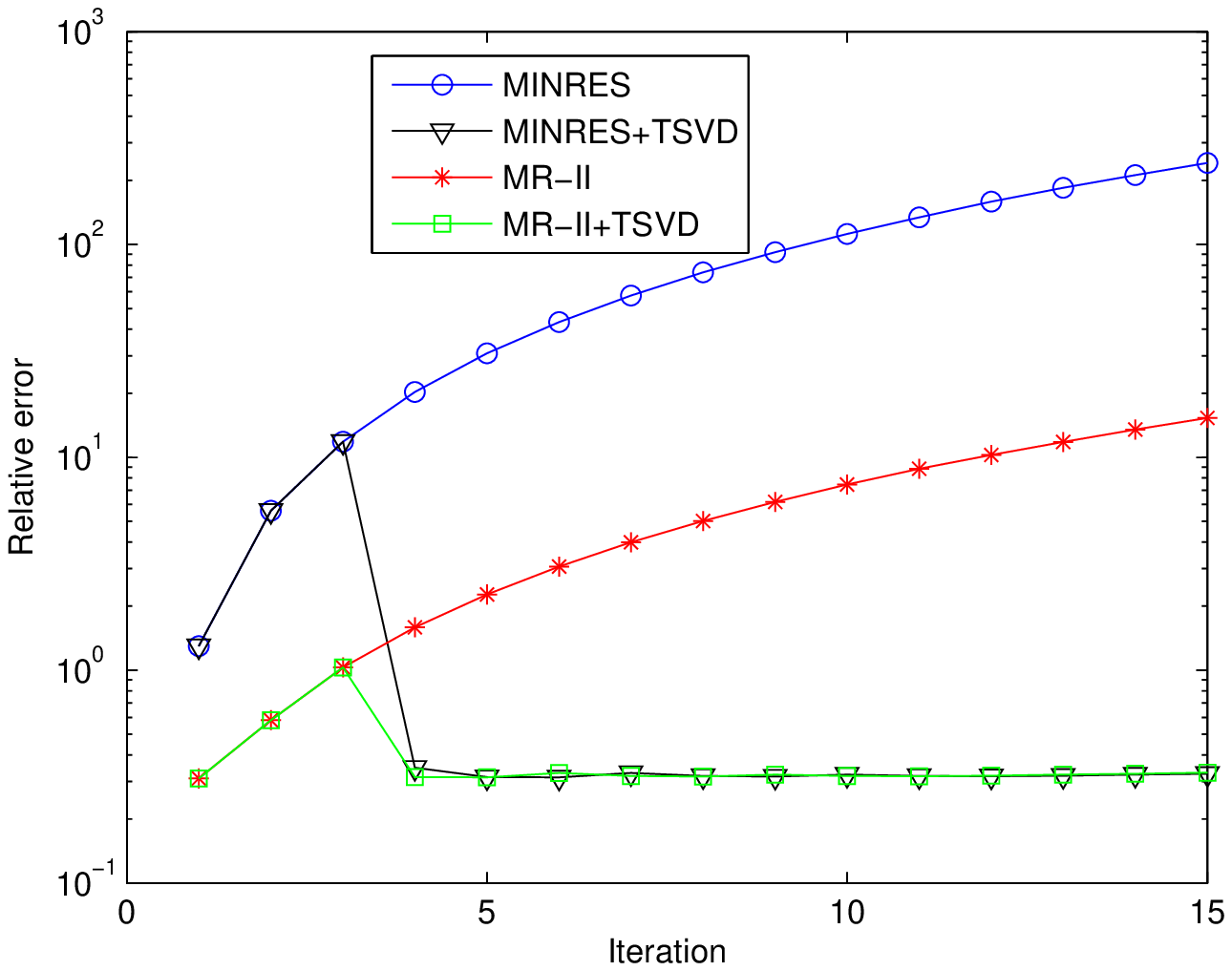}}
  \centerline{(a)}
\end{minipage}
\hfill
\begin{minipage}{0.48\linewidth}
  \centerline{\includegraphics[width=7.0cm,height=5cm]{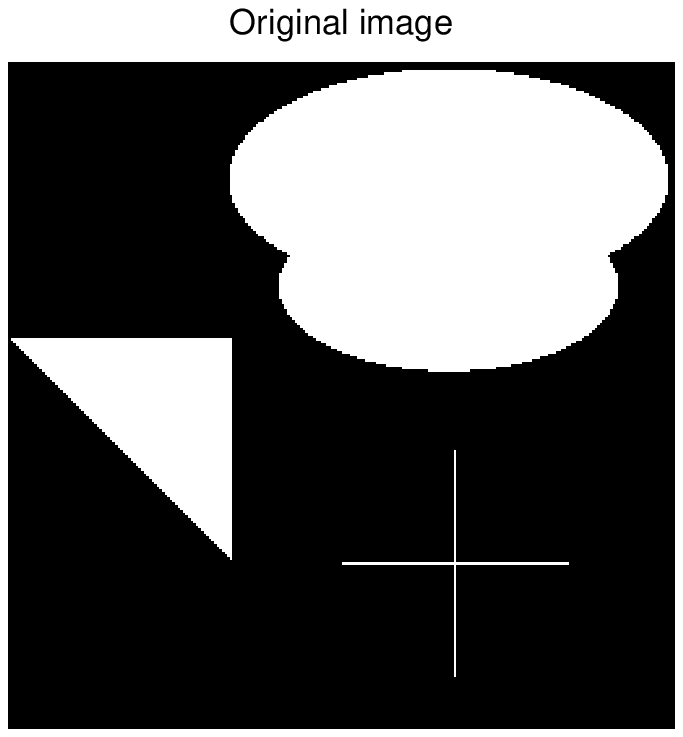}}
  \centerline{(b)}
\end{minipage}
\vfill
\begin{minipage}{0.48\linewidth}
  \centerline{\includegraphics[width=7.0cm,height=5cm]{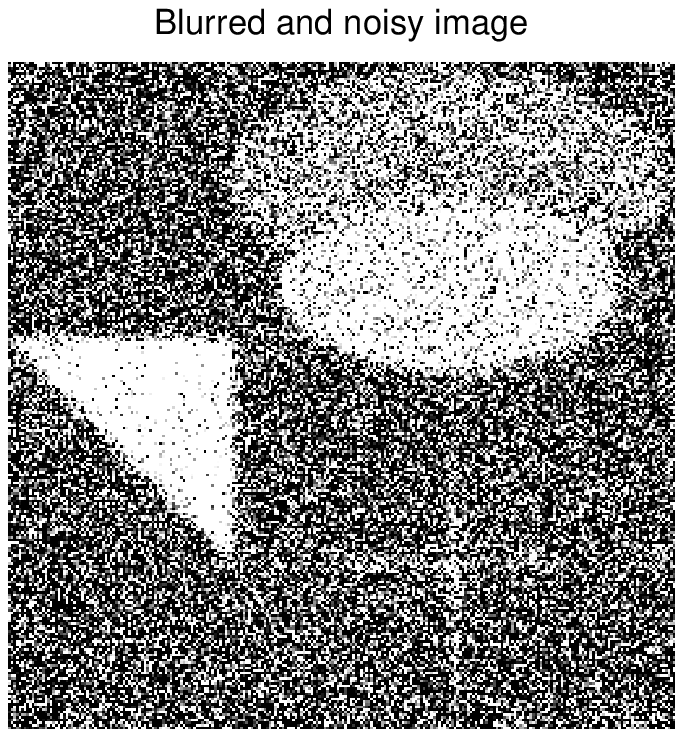}}
  \centerline{(c)}
\end{minipage}
\hfill
\begin{minipage}{0.48\linewidth}
  \centerline{\includegraphics[width=7.0cm,height=5cm]{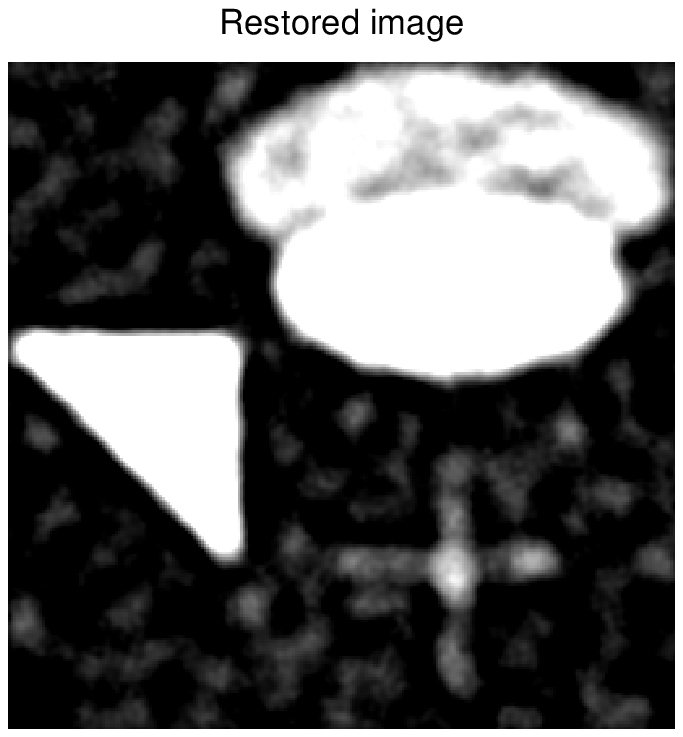}}
  \centerline{(d)}
\end{minipage}
\caption{(a): The relative errors $\|x^{(k)}-x_{true}\|/\|x_{true}\|$
with respect to MINRES, hybrid MINRES, MR-II,
and hybrid MR-II; (b): The original image; (c): The blurred and noisy image;
(d): The restored image with $\mathsf{band=7}, \mathsf{sigma=2}$.}
\label{fig12}
\end{figure}

Figure~\ref{fig11} (a) shows that MINRES and MR-II have the partial
regularization for the mildly ill-posed problem blur.
The semi-convergence of
the two methods occurs at the very first iteration, then regularized
solutions are progressively deteriorated,
while the hybrid MINRES finds the best possible regularized solution at
iteration $k=9$ and the hybrid MR-II does so at $k=8$.
Moreover, we see that the hybrid MINRES and MR-II reaches the same minimum
error level. Figure~\ref{fig11} exhibits the restoration performance,
where the restored image is
chosen by the regularized solution at the iteration where the
hybrid MR-II first reaches the minimum error level.
We observe from Figure~\ref{fig11} (d) that the outline of original
image is restored quite well by the restored image.

From Figure~\ref{fig12} (a), we see that the semi-convergence of
MR-II occurs at the first iteration and the regularized solution at
this iteration is as accurate as those obtained by the hybrid MR-II and
MINRES for the moderately ill-posed problem blur. Therefore, MR-II has
the full regularization for this problem. In contrast,
MINRES has only the partial regularization because its regularized
solution at semi-convergence is much less accurate than that obtained by MR-II.
In addition, we observe that the hybrid MINRES and the hybrid MR-II simply make
the regularized solutions almost stabilize with the minimum error.
Figure~\ref{fig12} (d) exhibits the restored image, which is a
good approximation to the original image.

\section{Conclusions}\label{SectionCon}

For large scale symmetric discrete linear ill-posed problems, MINRES and MR-II are
natural alternatives to LSQR. Our theory and
experiments have shown that MINRES has only the partial regularization and its
hybrid variant is needed to find best possible regularized solutions,
independent of the degree of ill-posedness. We have proved that MR-II has
better regularizing effects for severely and moderately ill-posed problems than
for mildly ill-posed problems, and it generally has only the partial
regularization for mildly ill-posed problems. We have shown that
although MR-II is a better regularization method than MINRES, the $k$th
MINRES iterate is always more accurate than the $(k-1)$th MR-II iterate until the
semi-convergence of MINRES. We have also established estimates for the
entries generated by the Lanczos process working on $\mathcal{K}(A,Ab)$, showing
how fast they decay. All these results have been
confirmed numerically. Remarkably,
stronger than our theory predicts, we have numerically demonstrated that MR-II
has the full regularization for severely and moderately ill-posed problems
and can compute best possible regularized solutions.
As a comparison of MR-II and
LSQR for a general symmetric ill-posed problem, our theory experiments have
indicated that two methods have very similar
regularizing effects but MR-II is twice as efficient as LSQR, so do their
hybrid variants. Therefore, for a large scale symmetric problem \eqref{eq1},
MR-II may be preferable to LSQR.

Some problems need to be further considered.
As we have seen, more appealing is a sharp estimate
for $\|\Delta_k\|$ other than $\|\Delta_k\|_F$.
The quantity $\|\sin\Theta(\mathcal{V}_k,\mathcal{V}_k^s)\|$
needs a more subtle analysis and plays a central role
in accurately estimating the accuracy $\gamma_k$ of the rank $k$
approximation generated by the Lanczos process working on $\mathcal{K}(A,Ab)$.
As we have elaborated, studying how near $\gamma_k$ is to $\sigma_{k+1}$ is a
central problem that completely understands the regularizing
effects of MR-II. Our bounds in Theorems~\ref{thm2}--\ref{thm3} are
less sharp and need to be improved on.



\end{document}